\documentclass[12pt,final]{siamltex}
\usepackage[dvips]{graphics}
\usepackage{amssymb}
\usepackage{doublespace}
\usepackage{times}

\makeatletter
\newcommand{\romref}[1]{%
 \expandafter\ifx\csname r@#1\endcsname\relax%
   \message{There were undefined references!!}%
   \textbf{!!}%
 \else%
   \romannumeral\ref{#1}%
 \fi%
}
\makeatother

\newcommand{\Fad}{\widehat{F}_{\mathrm{ad}}}

\newcommand{\Spec}{\mathop{\mathrm{Spec}}\nolimits}

\newcommand{\adj}{\mathop{\mathrm{adj}}\nolimits}
\newcommand{\Kf}{K\!_f}
\newcommand{\IR}{\mathbb{R}}

\newcommand{\bmatrix}[1]{\left[\matrix{#1}\right]}

\renewcommand{\Bigl}[1]{\mbox{\Large$\bigl#1$}}
\renewcommand{\Bigr}[1]{\mbox{\Large$\bigr#1$}}

\newtheorem{remarkaux}[theorem]{{\textit{Remark}}}
\newenvironment{remark}{\begin{remarkaux}\upshape}{\end{remarkaux}}
\newtheorem{assumptionaux}[theorem]{{\textit{Assumption}}}
\newenvironment{assumption}{\begin{assumptionaux}\upshape}{\end{assumptionaux}}
\newtheorem{exampleaux}[theorem]{{\textit{Example}}}
\newenvironment{example}{\begin{exampleaux}\upshape}{\end{exampleaux}}

\title{Feedback Stabilization over Commutative Rings:\\
       Further study of the coordinate-free approach}
\author{Kazuyoshi~Mori\footnotemark[1] \footnotemark[2] \and Kenichi~Abe\footnotemark[1]}

\begin{document}
\maketitle
\renewcommand{\thefootnote}{\fnsymbol{footnote}}
\footnotetext[1]{Department of Electrical Engineering, 
                 Faculty of Engineering, 
                 Tohoku University, Sendai 980-8579, JAPAN 
                ({\tt Kazuyoshi.MORI@IEEE.ORG, 
                      abe@abe.ecei.tohoku.ac.jp})}
\footnotetext[2]{This paper was partially written
while the first author was visiting the 
Institut de Recherche en Cybern\'etique de Nantes, UMR 6597, France
(September 1998 -- June 1999).}
\renewcommand{\thefootnote}{\arabic{footnote}}

  \centerline{\textsl{\today}}

\begin{abstract}
This paper is concerned with the coordinate-free approach to control
systems.  The coordinate-free approach is a~factorization approach but
does not require the coprime factorizations of the plant.  We present
two criteria for feedback stabilizability for MIMO systems in which
transfer functions belong to the total rings of fractions of
commutative rings.  Both of them are generalizations of Sule's results
in $[$SIAM J.\,Control Optim., \textbf{32}--6, 1675--1695(1994)$]$.
The first criterion is expressed in terms of modules generated from
a~causal plant and does not require the plant to be strictly causal.
It shows that if the plant is stabilizable, the modules are
projective.  The other criterion is expressed in terms of ideals
called generalized elementary factors.  This gives the stabilizability
of a~causal plant in terms of the coprimeness of the generalized
elementary factors.  As an example, a~discrete finite-time delay
system is considered.
\end{abstract}

\begin{keywords} 
  Linear systems, Feedback stabilization, 
  Coprime factorization over commutative rings
\end{keywords}

\begin{AMS}
93C05, %
93D15, %
93B50, %
93B25  %
\end{AMS}

  \makeatletter
  \begin{@abssec}{Abbreviated Title}
    Feedback Stabilization over Commutative Rings
  \end{@abssec}
  \makeatother

\pagestyle{myheadings}
\thispagestyle{plain}

  \markboth{\relax}{\relax}

\section{Introduction}\label{S:Introduction}
In this paper we are concerned with the coordinate-free approach to
control systems.  This approach is a~factorization approach but does
not require the coprime factorizations of the plant.

The factorization approach to control systems has the advantage that
it embraces, within a~single framework, numerous linear systems such
as continuous-time as well as discrete-time systems, lumped as well as
distributed systems, $1$-D as well as~$n$-D systems,
etc.\cite{bib:vidyasagar82a}.  This factorization approach was
patterned after Desoer \emph{et al.}\cite{bib:desoer80a} and
Vidyasagar \emph{et al.}\cite{bib:vidyasagar82a}.  In this approach,
when problems such as feedback stabilization are studied, one can
{focus} on the key aspects of the problem under study rather than be
distracted by the special features of a~particular class of linear
systems.  A~transfer function of this approach is given as the ratio
of two stable causal transfer functions and the set of stable causal
transfer functions forms a~commutative ring.  For a~long time, the
theory of the factorization approach had been founded on the coprime
factorizability of transfer matrices, which is satisfied in the case
where the set of stable causal transfer functions is such
a~commutative ring as a~Euclidean domain, a~principal ideal, or a~B\'ezout domain.  However, Anantharam in~\cite{bib:anantharam85a}
showed that there exist models in which some stabilizable plants do
not have right-/left-coprime factorizations.

Recently, Shankar and Sule in~\cite{bib:shankar92a} have presented
a~theory of feedback stabilization for single-input single-output
(SISO) transfer functions having fractions over general integral
domains.  Moreover, Sule in~\cite{bib:sule94a,bib:sule98a} has
presented a~theory of the feedback stabilization of strictly causal
plants for multi-input multi-output (MIMO) transfer matrices, in which
transfer functions belong to the total rings of fractions of
commutative rings, with some restrictions.  Their approach to the
control systems is called a ``coordinate-free
approach''(\cite[p.15]{bib:shankar92a}) in the sense that they do not
require the coprime factorizability of transfer matrices.

The main contribution of this paper consists of providing two criteria
for feedback stabilizability for MIMO systems in which transfer
functions belong to the total rings of fractions of commutative rings:
the first criterion is expressed in terms of modules
((\romref{Th:3.3:ii}) of~Theorem\,\ref{Th:3.3}) and the other in terms
of ideals called generalized elementary factors ((\romref{Th:3.3:iii})
of~Theorem\,\ref{Th:3.3}).  They are more general than Sule's results
in the following sense:~(i) our results do not require that plants are
strictly causal;~(ii) we do not employ the restriction of commutative
rings.
Further, we will not use the theory of algebraic geometry.

The paper is organized as follows.  In~\S\,\ref{S:Preliminary}, we
give mathematical preliminaries, set up the feedback stabilization
problem, present the previous results, and define the causality of the
transfer functions.  In~\S\,\ref{S:MainResults}, we state our main
results.  As a~preface to our main results, we also introduce there
the notion of the generalized elementary factor of~a plant.
In~\S\,\ref{S:IR}, we give intermediate results which we will utilize
in the proof of the main theorem.  In~\S\,\ref{S:ProofMainResults}, we
prove our main theorem.  In~\S\,\ref{S:CSC}, we discuss the causality
of the stabilizing controllers.  Also, in order to make the contents
clear, we present examples concerning a~discrete finite-time delay
system in~\S\,\ref{S:MainResults}, \ref{S:IR},
\ref{S:ProofMainResults} in series.

\section{Preliminaries}
\label{S:Preliminary}
In the following we begin by introducing the notations of commutative
rings, matrices, and modules, commonly used in this paper.  Then we
give in order the formulation of the feedback stabilization problem,
the previous results, and the causality of transfer functions.

\subsection{Notations}\label{SS:Notations}
\paragraph{Commutative Rings}
In this paper, we consider that any commutative ring has the
identity~$1$ different from zero.  Let~${\cal R}$ denote a~commutative ring.
A \emph{zerodivisor} in~${\cal R}$ is an element~$x$ for which there exists
a~nonzero~$y$ such that $xy=0$.  In particular, a~zerodivisor~$x$ is
said to be \emph{nilpotent}
if $x^n=0$ for some positive integer~$n$.  The set of all nilpotent
elements in~${\cal R}$, which is an ideal, is called the \emph{nilradical}
of~${\cal R}$.  A \emph{nonzerodivisor} in~${\cal R}$ is an element which is not
a~zerodivisor.  The total ring of fractions of~${\cal R}$ is denoted by~${\cal F}
({\cal R})$.

The set of all prime ideals of~${\cal R}$ is called the \emph{prime
spectrum} of~${\cal R}$ and is denoted by $\Spec{\cal R}$.  The prime spectrum
of~${\cal R}$ is said to be \emph{irreducible} as a~topological space if
every non-empty open set is dense in $\Spec{\cal R}$.

We will consider that the set of stable causal transfer functions is
a~commutative ring, denoted by~${\cal A}$.  Further, we will use the
following three kinds of rings of fractions:
\begin{romannum}
\item
The first one appears as the total ring of fractions of~${\cal A}$, which is
denoted by~${\cal F}({\cal A})$ or simply by~${\cal F}$; that is, ${\cal F}=\{ n/d\,|\, n,d\in{\cal A}$,
$d$ is a~nonzerodivisor$\}$.  This will be considered to be the set of
all possible transfer functions.  If~the commutative ring~${\cal A}$ is an
integral domain,~${\cal F}$ becomes a~field of fractions of~${\cal A}$.  However,
if~${\cal A}$ is not an integral domain, then~${\cal F}$ is not a~field, because
any zerodivisor of~${\cal F}$ is not a~unit.
\item
The second one is associated with the set of powers of a~nonzero
element of~${\cal A}$.  Suppose that~$f$ denotes a~nonzero element of~${\cal A}$.
Given a~set $S_f=\{ 1,f,f^2,\ldots\}$, which is a~multiplicative subset
of~${\cal A}$, we denote by~${\cal A}_f$ the ring of fractions of~${\cal A}$ with
respect to the multiplicative subset~$S_f$; that is, ${\cal A}_f=\{ n/d\,|\, n\in
{\cal A},~d\in S_f\}$.  We point out two facts here: (a) In the case where~$f$
is nilpotent, ${\cal A}_f$ becomes isomorphic to $\{ 0\}$. 
(b) In the case where~$f$ is a~zerodivisor, even if the equality
$a/1=b/1$ holds over~${\cal A}_f$ with $a,b\in{\cal A}$, we cannot say in general
that $a=b$ over~${\cal A}$; alternatively, $a=b+z$ over~${\cal A}$ holds with some
zerodivisor~$z$ of~${\cal A}$ such that $zf^{\omega}=0$ with a~sufficiently large
integer~$\omega$.
\item
The last one is the total ring of fractions of~${\cal A}_f$, which is
denoted by ${\cal F}({\cal A}_f)$; that is, ${\cal F}({\cal A}_f)=\{ n/d\,|\, n,d\in{\cal A}_f$, $d$ is
a~nonzerodivisor of ${\cal A}_f\}$.  If~$f$ is a~nonzerodivisor of~${\cal A}$, ${\cal F}({\cal A}_f) $
coincides with the total ring of fractions of~${\cal A}$.  Otherwise, they
may not coincide.
\end{romannum}

The reader is referred to~Chapter\,3 of~\cite{bib:atiyah69a} for the ring
of fractions and to~Chapter\,1 of~\cite{bib:atiyah69a} for the prime
spectrum.

In the rest of the paper, we will use~${\cal R}$ as an unspecified commutative
ring and mainly suppose that~${\cal R}$ denotes either~${\cal A}$ or~${\cal A}_f$.

\paragraph{Matrices}
Suppose that~$x$ and~$y$ denote sizes of matrices.

The set of matrices over~${\cal R}$ of size $x\times y$ is denoted by ${\cal R}^{x\times
y}$.  In particular, the set of square matrices over~${\cal R}$ of size~$x$
is denoted by~$({\cal R})_x$.  A~square matrix is called \emph{singular}
over~${\cal R}$ if its determinant is a~zerodivisor of~${\cal R}$, and
\emph{nonsingular} otherwise. 
The identity and the zero matrices are denoted by~$E_x$ and $O_{x\times
y}$, respectively, if the sizes are required, otherwise they are
denoted simply by~$E$ and~$O$.
For a~matrix~$A$ over~${\cal R}$, the inverse matrix of~$A$ is denoted
by~$A^{-1}$ provided that $\det(A)$ is a~unit of ${\cal F}({\cal R})$.  The ideal
generated by~${\cal R}$-linear combination of all minors of size~$m$ of
a~matrix~$A$ is denoted by $I_{m{\cal R}}(A)$.

Matrices~$A$ and~$B$ over~${\cal R}$ are \emph{right-coprime over~${\cal R}$} if
there exist matrices~$\widetilde{X}$ and~$\widetilde{Y}$ over~${\cal R}$
such that
  $\widetilde{X}A+\widetilde{Y}B=E$. 
Analogously, matrices~$\widetilde{A}$ and~$\widetilde{B}$ over~${\cal R}$
are \emph{left-coprime over~${\cal R}$} if there exist matrices~$X$ and~$Y$
over~${\cal R}$ such that
  $\widetilde{A}X+\widetilde{B}Y=E$. 
Note that, in the sense of the above definition, even if two matrices
have no common right-$\Bigl($left-$\Bigr)$divisors except invertible
matrices, they may not be right-$\Bigl($left-$\Bigr)$coprime over~${\cal R}
$.  (For example, two matrices $\bmatrix{z_1}$ and $\bmatrix{z_2}$ of
size~$1\times 1$ over the bivariate polynomial ring $\mathbb{R}[z_1,z_2]$
over the real field $\mathbb{R}$ are neither right- nor left-coprime
over $\mathbb{R}[z_1,z_2]$ in our setting.)  Further, a~pair $(N,D)$
of matrices~$N$ and~$D$ is said to be a \emph{right-coprime
factorization of~$P$ over~${\cal R}$} if~(i) the matrix~$D$ is nonsingular
over~${\cal R}$,~(ii) $P=ND^{-1}$ over ${\cal F}({\cal R})$, and~(iii)~$N$ and~$D$ are
right-coprime over~${\cal R}$.  Also, a~pair $(\widetilde{N},\widetilde{D})$
of matrices $\widetilde{N}$ and $\widetilde{D}$ is said to be a
\emph{left-coprime factorization of~$P$ over~${\cal R}$}
if~(i)~$\widetilde{D}$ is nonsingular over~${\cal R}$,~(ii)
$P=\widetilde{D}^{-1}\widetilde{N}$ over ${\cal F}({\cal R})$,
and~(iii)~$\widetilde{N}$ and $\widetilde{D}$ are left-coprime over~$
{\cal R}$.  As we have seen, in the case where a~matrix is potentially used
to express \emph{left} fractional form and/or \emph{left} coprimeness,
we usually attach a~tilde `$\widetilde{~~}$' to a~symbol; for example
$\widetilde{N}$, $\widetilde{D}$ for
$P=\widetilde{D}^{-1}\widetilde{N}$ and $\widetilde{Y}$,
$\widetilde{X}$ for $\widetilde{Y}N+\widetilde{X}D=E$.

\paragraph{Modules}
For~a matrix~$A$ over~${\cal R}$, we denote by $M_r(A)$
\,$\Bigl(M_c(A)\Bigr)$ the~${\cal R}$-module generated by rows
\,$\Bigl($columns$\Bigr)$ of~$A$.

Suppose that~$A$, $B$, $\widetilde{A}$, $\widetilde{B}$ are matrices
over~${\cal R}$ and~$X$ is a~matrix over ${\cal F}({\cal R})$ such that
$X=AB^{-1}=\widetilde{B}^{-1}\widetilde{A}$ with~$B$ and
$\widetilde{B}$ being nonsingular.  Then the ${\cal R}$-module
$M_r(\bmatrix{A^t & B^t}^t)$ $\Bigl(M_c(\bmatrix{\widetilde{A} &
\widetilde{B}})\Bigr)$ is uniquely determined up to isomorphism with
respect to any choice of fractions $AB^{-1}$
$\Bigl(\widetilde{B}^{-1}\widetilde{A}\Bigr)$ of~$X$ as shown in
Lemma\,\ref{L:2.1} below.  Thus for a~matrix~$X$ over~${\cal F}({\cal R}) $, we
denote by~${\cal T}_{X,{\cal R}} $ and~${\cal W}_{X,{\cal R}}$ the modules $M_r(\bmatrix{A^t &
B^t}^t)$ and $M_c(\bmatrix{\widetilde{A} & \widetilde{B}})$,
respectively.  If ${\cal R}={\cal A}$, we write simply~${\cal T}_X$ and~${\cal W}_X$ for~${\cal T}_
{X,{\cal A}}$ and~${\cal W}_{X,{\cal A}} $, respectively.  We will use, for example,
the notations~${\cal T}_P$, ${\cal W}_ P$, ${\cal T}_C$, and~${\cal W}_C$ for the matrices~$P$
and~$C$ over~${\cal F}$.

An~${\cal R}$-module~$M$ is called \emph{free} if it has a~basis, that is,
a~linearly independent system of generators.  The \emph{rank} of
a~free~${\cal R}$-module~$M$ is equal to the cardinality of a~basis of~$M$,
which is independent of the basis chosen.  An ${\cal R}$-module~$M$ is
called \emph{projective} if it is a~direct summand of a~free~${\cal R}
$-module, that is, there is a~module~$N$ such that $M\oplus N$ is free.
The reader is referred to~Chapter\,2 of~\cite{bib:atiyah69a} for the
module theory.

\begin{lemma}\label{L:2.1}
Suppose that~$X$ is a~matrix over ${\cal F}({\cal R})$ and is expressed in the
form of a~fraction $X=AB^{-1}=\widetilde{B}^{-1} \widetilde{A}$ with
some matrices~$A$, $B$, $\widetilde{A}$, $\widetilde{B}$ over~${\cal R}$.
Then the ${\cal R}$-module $M_r(\bmatrix{A^t & B^t}^t)$
$\Bigl(M_c(\bmatrix{\widetilde{A} & \widetilde{B}})\Bigr)$ is uniquely
determined up to isomorphism with respect to any choice of fractions
$AB^{-1}$ $\Bigl(\widetilde{B}^{-1}\widetilde{A}\Bigr)$ of~$X$.
\end{lemma}
\begin{proof}
Without loss of generality, it is sufficient to show that
$M_r(\bmatrix{A_1^t & b_1E}^t) \simeq M_r(\bmatrix{A_2^t & B_2^t}^t)$,
where $b_1\in{\cal R}$ and $A_1(b_1E)^{-1}=A_2B_2^{-1}$.  Since~$b_1$ is
a~nonzerodivisor and~$B_2$ is nonsingular, we have $
M_r(\bmatrix{A_1^t & b_1E}^t) \simeq M_r(\bmatrix{A_1^t & b_1E}^tB_2) \simeq
M_r(\bmatrix{A_2^t & B_2^t}^tb_1E)\simeq M_r(\bmatrix{A_2^t & B_2^t}^t)$.
The other isomorphism can be proved analogously.  \qquad
\end{proof}

\subsection{Feedback Stabilization Problem}\label{SS:FSProblem}
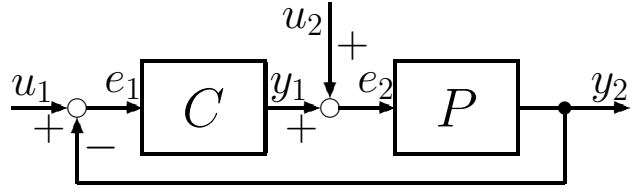
\begin{figure}
\begin{center}
%
\setlength{\unitlength}{0.08mm}
\begin{picture}(1160.000000,400.000000)(-50.000000,-50.000000)
\linethickness{0.4mm}
\put(220.000000,50.000000){\framebox(200.000000,150.000000){\LARGE$C$}}
\put(640.000000,50.000000){\framebox(200.000000,150.000000){\LARGE$P$}}
\put(110.000000,125.000000){\circle{30.000000}}
\put(530.000000,125.000000){\circle{30.000000}}
\put(0.000000,125.000000){\vector(1,0){85.000000}}
\put(85.000000,129.000000){\vector(1,0){0}}
\put(85.000000,121.000000){\vector(1,0){0}}
\put(95.000000,125.000000){\vector(1,0){0}}
\put(125.000000,125.000000){\vector(1,0){85.000000}}
\put(210.000000,129.000000){\vector(1,0){0}}
\put(210.000000,121.000000){\vector(1,0){0}}
\put(220.000000,125.000000){\vector(1,0){0}}
\put(420.000000,125.000000){\vector(1,0){85.000000}}
\put(505.000000,129.000000){\vector(1,0){0}}
\put(505.000000,121.000000){\vector(1,0){0}}
\put(515.000000,125.000000){\vector(1,0){0}}
\put(545.000000,125.000000){\vector(1,0){85.000000}}
\put(630.000000,129.000000){\vector(1,0){0}}
\put(630.000000,121.000000){\vector(1,0){0}}
\put(640.000000,125.000000){\vector(1,0){0}}
\put(840.000000,125.000000){\vector(1,0){180.000000}}
\put(1020.000000,129.000000){\vector(1,0){0}}
\put(1020.000000,121.000000){\vector(1,0){0}}
\put(1030.000000,125.000000){\vector(1,0){0}}
\put(110.000000,0.000000){\vector(0,1){100.000000}}
\put(106.000000,100.000000){\vector(0,1){0}}
\put(114.000000,100.000000){\vector(0,1){0}}
\put(110.000000,110.000000){\vector(0,1){0}}
\put(110.000000,0.000000){\line(1,0){810.000000}}
\put(920.000000,0.000000){\line(0,1){125.000000}}
\put(530.000000,300.000000){\vector(0,-1){150.000000}}
\put(526.000000,150.000000){\vector(0,-1){0}}
\put(534.000000,150.000000){\vector(0,-1){0}}
\put(530.000000,140.000000){\vector(0,-1){0}}
\put(0.000000,135.000000){\makebox(0,0)[lb]{\Large$u_1$}}
\put(218.000000,140.000000){\makebox(0,0)[rb]{\Large$e_1$}}
\put(90.000000,115.000000){\makebox(0,0)[rt]{\Large$+$}}
\put(120.000000,85.000000){\makebox(0,0)[lt]{\Large$-$}}
\put(430.000000,135.000000){\makebox(0,0)[lb]{\Large$y_1$}}
\put(520.000000,290.000000){\makebox(0,0)[rt]{\Large$u_2$}}
\put(638.000000,140.000000){\makebox(0,0)[rb]{\Large$e_2$}}
\put(1028.000000,135.000000){\makebox(0,0)[rb]{\Large$y_2$}}
\put(510.000000,115.000000){\makebox(0,0)[rt]{\Large$+$}}
\put(540.000000,205.000000){\makebox(0,0)[lb]{\Large$+$}}
\put(920.000000,125.000000){\circle*{22.500000}}
\end{picture}

%
  \caption{Feedback system~$\Sigma$.}\label{Fig:FeedbackSystem}
\end{center}
\end{figure}

The stabilization problem considered in this paper follows that of
Sule in~\cite{bib:sule94a} who considers the feedback system~$\Sigma
$~\cite[Ch.5, Figure\,5.1]{bib:vidyasagar85a} as in
Figure\,\ref{Fig:FeedbackSystem}.
For further details the reader is referred
to~\cite{bib:vidyasagar85a}.  Let a~commutative ring~${\cal A}$ represent
the set of \emph{stable causal} transfer functions.  The total ring of
fractions of~${\cal A}$, denoted by~${\cal F}$, consists of \emph{all} possible
transfer functions.  The set of matrices of size $x\times y$ over~${\cal A}$,
denoted by ${\cal A}^{x\times y}$, coincides with the set of stable causal
transfer matrices of size $x\times y$.  Also the set of matrices of size $x
\times y$ over~${\cal F}$, denoted by ${\cal F}^{x\times y}$, coincides with all possible
transfer matrices of size $x\times y$.  Throughout the paper, the plant we
consider has~$m$ inputs and~$n$ outputs, and its transfer matrix,
which itself is also called simply a \emph{plant}, is denoted by~$P$
and belongs to~${\cal F}^{n\times m}$.  We will occasionally consider matrices
over~${\cal A}$ $\Bigl({\cal F}\Bigr)$ as ones over~${\cal A}_f$ or~${\cal F}$ $\Bigl({\cal F}({\cal A}_
f)\Bigr)$ by natural mapping.

\begin{definition}\label{D:Unreferred01}
Define $\Fad$ by
\[
  \Fad=\{(X,Y)\in{\cal F}^{x\times y}\times{\cal F}^{y\times x}\,|\,
            \begin{array}[t]{l}
               \det(E_x+XY)\mbox{ is a~unit of }{\cal F}, \\
              \mbox{$x$ and~$y$ are positive integers}\}.
            \end{array}
\]
For $P\in{\cal F}^{n\times m}$ and $C\in{\cal F}^{m\times n}$, the matrix $H(P,C)\in({\cal F})_{m+n}$
is defined by
\begin{equation}\label{E:H(P,C)}
  H(P,C) =
   \bmatrix{
    (E_n+PC)^{-1}  &  -P(E_m+CP)^{-1} \cr
    C(E_n+PC)^{-1} & (E_m+CP)^{-1}}
\end{equation}
provided $(P,C)\in\Fad$.  This $H(P,C)$ is the transfer matrix from
$\bmatrix{u_1^t & u_2^t}^t$ to $\bmatrix{e_1^t & e_2^t}^t$ of the
feedback system~$\Sigma$.  If~(i) $(P,C)\in\Fad$ and~(ii) $H(P,C)\in({\cal A})_
{m+n}$, then we say that the plant~$P$ is \emph{stabilizable}, $P$ is
\emph{stabilized} by~$C$, and~$C$ is a \emph{stabilizing controller}
of~$P$.
\end{definition}

\subsection{Previous Results}
\label{SS:PrevRes}
In~\cite{bib:sule94a} Sule gave the results of the feedback
stabilizability.  We show them after introducing the notion of the
elementary factor which is used to state his results.
\begin{definition}\label{D:EF}
\textup{(Elementary Factors \cite[p.1689]{bib:sule94a})~} Assume
that~${\cal A}$ is a~unique factorization domain.  Denote by~$T$ the matrix
$\bmatrix{N^t&dE_m}^t$ and by~$W$ the matrix $\bmatrix{N & dE_n}$
over~${\cal A}$, where $P=Nd^{-1}$ with $N\in{\cal A}^{n\times m}$, $d\in{\cal A}$.  Let $
\{ T_1,T_2,\ldots, T_r\}$ be the family of all nonsingular $m\times m$
submatrices of the matrix~$T$, and for each index~$i$, let~$f_i$ be
the radical of the least common multiple of all the denominators of
$TT_i^{-1}$.  The family $F=\{ f_1,f_2,\ldots, f_r\}$ is called the family
of \emph{elementary factors of the matrix $T$}.  Analogously let $
\{ W_1,W_2,\ldots, W_r\}$ be the family of all nonsingular $n\times n$
submatrices of the matrix~$W$, and for each index~$j$, let~$g_j$ be
the radical of the least common multiple of all the denominators of
$W_j^{-1}W$.  Let $G=\{ g_1,g_2,\ldots, g_l\}$ denote the family of
\emph{elementary factors of the transposed matrix $W^t$}.  Now let $H=
\{ f_ig_j\,|\, i=1,\ldots, r,~j=1,\ldots, l\}$.  This family~$H$ is called the
\emph{elementary factor of the transfer matrix $P$}.
\end{definition}

Then, Sule's two elegant results can be rewritten as follows.  The
first result assumes that the prime spectrum of~${\cal A}$ is irreducible.
The second one assumes that~${\cal A}$ is a~unique factorization domain.
\begin{theorem}\label{Sule94:Th:1}
\textup{(Theorem\,1 of~\cite{bib:sule94a})~} Suppose that the prime
spectrum of~${\cal A}$ is irreducible.  Further suppose that a~plant~$P$ of
size $n\times m$ is strictly causal, where the notion of the strictly
causal is defined as in~\cite{bib:sule98a} (rather
than~\cite{bib:sule94a}).  Then the plant~$P$ is stabilizable if and
only if the following conditions are satisfied:
\begin{romannum}
\item
  The module~${\cal T}_P$ is projective of rank~$m$.
\item
  The module~${\cal W}_P$ is projective of rank~$n$.
\qquad\endproof
\end{romannum}
\end{theorem}

\noindent
Recall that for a~matrix~$X$ over~${\cal F}$, we use the notations~${\cal T}_X$
and~${\cal W}_X$ to denote~${\cal A}$-modules generated by using the matrix~$X$.
Further it should be noted that the definitions of ${\cal T}_P$ and ${\cal W}_P$
in~\cite{bib:sule94a} are slightly different from those of this paper.
Nevertheless this is not a problem by virtue of Lemma\,\ref{L:2.1}.
\begin{theorem}\label{Sule94:Th:4}
\textup{(Theorem\,4 of~\cite{bib:sule94a})~} Suppose that~${\cal A}$ is
a~unique factorization domain.  The plant~$P$ is stabilizable if and
only if the elementary factors of~$P$ are coprime, that is,
$\sum_{h\in H}(h)={\cal A}$.
\qquad\endproof
\end{theorem}

\subsection{Causality of Given Plants}
Here we define the causality of transfer functions, which is an
important physical constraint, used in this paper.  We employ the
definition of causality from Vidyasagar~\emph{et
al.}\cite[Definition\,3.1]{bib:vidyasagar82a} and introduce two
terminologies used later frequently.
\begin{definition}\label{D:Causal}
Let~${\cal Z}$ be a~prime ideal of~${\cal A}$, with ${\cal Z}\neq{\cal A}$, including all
zerodivisors.  Define the subsets~${\cal P}$ and ${\cal P}_{\textrm{s}}$ of~${\cal F}$
as follows:
\[
  {\cal P}=\{ n/d\in{\cal F}\,|\, n\in{\cal A},~ d\in{\cal A}\backslash{\cal Z}\},{\ \ }
  {\cal P}_{\textrm{s}}=\{ n/d\in{\cal F}\,|\, n\in{\cal Z},~ d\in{\cal A}\backslash{\cal Z}\}.
\]
A~transfer function in~${\cal P}$ $\Bigl({\cal P}_{\textrm{s}}\Bigr)$ is called
\emph{causal} $\Bigl($\emph{strictly causal}$\Bigr)$.  Similarly, if
every entry of a~transfer matrix over~${\cal F}$ is in~${\cal P}$ $\Bigl({\cal P}_
{\textrm{s}}\Bigr)$, the transfer matrix is called \emph{causal}
$\Bigl($\emph{strictly causal}$\Bigr)$.  A~transfer matrix is said to
be \emph{${\cal Z}$-nonsingular} if the determinant is in ${\cal A}\backslash{\cal Z}$, and
\emph{${\cal Z}$-singular} otherwise.
\end{definition}

In~\cite{bib:vidyasagar82a}, the ideal~${\cal Z}$ is not restricted to
a~prime ideal in general.  On the other hand, in~\cite{bib:sule94a},
the set of the denominators of causal transfer functions is
a~multiplicatively closed subset of~${\cal A}$.  This property is natural
since the multiplication of two causal transfer functions should be
considered as causal one.  Note that this multiplicativity is
equivalent to~${\cal Z}$ being prime provided that ${\cal Z}$ is an ideal.
By following the multiplicativity, we consider in this paper that~${\cal Z}$
is prime.

In this paper, we do not assume that the prime spectrum of~${\cal A}$ is
irreducible and the plant~$P$ is strictly causal as
in~\cite{bib:sule94a}.  Alternatively, in the rest of the paper we
assume only the following:
\begin{assumption}\label{A:domainassumption}
The given plant is causal in the sense of Definition\,\ref{D:Causal}.
\qquad\endproof
\end{assumption}

One can represent a~causal plant~$P$ in the form of fractions
$P=ND^{-1}=\widetilde{D}^{-1}\widetilde{N}$, where the matrices~$N$,
$D$, $\widetilde{N}$, $\widetilde{D}$ are over~${\cal A}$, and the
matrices~$D$, $\widetilde{D}$ are ${\cal Z}$-nonsingular.

\section{Main Results}\label{S:MainResults}
To state our results precisely we define the notion of generalized
elementary factors, which is a~generalization of the elementary
factors in Definition\,\ref{D:EF}.  Then the main theorem will be
presented.

\paragraph{Generalized Elementary Factors}
Originally, the elementary factors have been defined over unique
factorization domains as in Definition\,\ref{D:EF}.  The authors have
enlarged this concept for integral domains\cite{bib:mori96a} and have
presented a~criterion for feedback stabilizability over integral
domains.  We enlarge this concept again in the case of commutative
rings.

Before stating the definition, we introduce several symbols used in
the definition and widely in the rest of this paper.  The symbol~${\cal I}$
denotes the set of all sets of~$m$ distinct integers between~$1$ and
$m+n$ (recall that~$m$ and~$n$ are the numbers of the inputs and the
outputs, respectively).  Normally, an element of~${\cal I}$ will be denoted
by~$I$, possibly with suffixes such as integers.  We will use an
element of~${\cal I}$ as a~suffix as well as a~set.  For~$I\in{\cal I}$, if
$i_1,\ldots, i_m$ are elements of~$I$ in ascending order, that is,
$i_a<i_b$ if $a<b$, then the symbol~$\Delta_I$ denotes the matrix whose
$(k,i_k)$-entry is~$1$ for $i_k\in I$ and zero otherwise.

\begin{definition}\label{D:GEF}
\textup{(Generalized Elementary Factors)~} Let $P\in{\cal F}^{n\times m}$, and~$N$
and~$D$ are matrices over~${\cal A}$ with $P=ND^{-1}$.  Denote by~$T$ the
matrix $\bmatrix{N^t&D^t}^t$.  For each $I\in{\cal I}$, define the ideal $\Lambda_
{P\!I}$ of~${\cal A}$ by
\[
  \Lambda_{P\!I}=\{\lambda\in{\cal A}\,|\,\exists K\in{\cal A}^{(m+n)\times m}\  \lambda T  =K \Delta_I T\}.
\]
We call the ideal $\Lambda_{P\!I}$ the \emph{generalized elementary factor}
of the plant~$P$ with respect to~$I$. 
\end{definition}

Whenever we use the symbol~$\Lambda$ with some suffix, it will denote
a~generalized elementary factor.  We will also frequently use the
symbols~$\lambda$ and~$\lambda_I$ with $I\in{\cal I}$ as particular elements of $\Lambda_
{P\!I}$.  Note that in Definitions\,\ref{D:EF} and~\ref{D:GEF}, the
matrices represented by~$T$ are different in general.  However this
difference is not a~problem since the generalized elementary factors
are independent of the choice of the fractions $ND^{-1}$ as shown
below.
\begin{lemma}\label{L:3.2}
For any $I\in{\cal I}$, the generalized elementary factor of the plant~$P$ with
respect to~$I$ is independent of the choice of matrices~$N$ and~$D$ over~$
{\cal A}$ satisfying $P=ND^{-1}$.
\end{lemma}
\begin{proof}
Let~$N$, $N'$, $D$ be matrices over~${\cal A}$ and~$d'$ be a~scalar of~${\cal A}$
such that $P=ND^{-1}=N'{d'}^{-1}$ hold.  Further, let
\begin{eqnarray*}
  \Lambda_{P\!I1}&=&\{\lambda\in{\cal A}\,|\,\exists K\in{\cal A}^{(m+n)\times m}\  
              \lambda\bmatrix{N^t&D^t}^t  =K \Delta_I \bmatrix{N^t&D^t}^t\},\\
  \Lambda_{P\!I2}&=&\{\lambda\in{\cal A}\,|\,\exists K\in{\cal A}^{(m+n)\times m}\  
              \lambda\bmatrix{{N'}^t&d'E_m}^t  =K \Delta_I \bmatrix{{N'}^t&d'E_m}^t\}.
\end{eqnarray*}
In order to prove this lemma it is sufficient to show that the ideals
$\Lambda_{P\!I1}$ and $\Lambda_{P\!I2}$ are equal.  Suppose that~$\lambda$ is an
element of $\Lambda_{P\!I1}$.  Then there exists a~matrix~$K$ such that
  $\lambda\bmatrix{N^t&D^t}^t =K \Delta_I \bmatrix{N^t&D^t}^t$.  
Multiplying $d'E_m$ on the right of both sides, we have  
  $\lambda\bmatrix{{N'}^t& d'E_m}^t D=$ $K \Delta_I \bmatrix{{N'}^t& d'E_m}^t D$.  
Since the matrix~$D$ is nonsingular, we have
  $\lambda\bmatrix{{N'}^t&d'E_m}^t =K \Delta_I \bmatrix{N^t&d'E_m}^t$,  
so that $\lambda\in \Lambda_{P\!I2}$, which means that $\Lambda_{P\!I1}\subset \Lambda_{P\!I2}$.  The
opposite inclusion relation $\Lambda_{P\!I1}\supset \Lambda_{P\!I2}$ can be proved
analogously.
\qquad
\end{proof}

Note also that for every~$I$ in~${\cal I}$, the generalized elementary
factor of the plant with respect to~$I$ is not empty since it contains
at least zero.  In the case where the set of stable causal transfer
functions is a~unique factorization domain, the generalized elementary
factor of the plant with the matrix $\Delta_IT$ being nonsingular becomes
a~principal ideal and the generator of its radical an elementary
factor of the matrix~$T$ (in Definition\,\ref{D:EF}) up to a~unit
multiple.  Analogously, the elementary factor of the matrix~$W$ (in
Definition\,\ref{D:EF}) corresponds to the generalized elementary
factor of the transposed plant~$P^t$.

\paragraph{Main Results}
We are now in position to state our main results.
\begin{theorem}\label{Th:3.3}
Consider a~causal plant~$P$.  Then the following statements are
equivalent:
\begin{romannum}
\item\label{Th:3.3:i}
The plant~$P$ is stabilizable.
\item\label{Th:3.3:ii} 
The ${\cal A}$-modules~${\cal T}_P$ and~${\cal W}_P$ are projective.
\item\label{Th:3.3:iii}
The set of all generalized elementary factors of~$P$ generates~${\cal A}$;
that is, 
\begin{equation}\label{E:Th:3.2:1}
  \sum_{I\in{\cal I}} \Lambda_{P\!I}={\cal A}.
\end{equation}
\qquad\endproof
\end{romannum}
\end{theorem}

In the theorem, (\romref{Th:3.3:ii}) and~(\romref{Th:3.3:iii}) are
criteria for feedback stabilizability.  Comparing the theorem above
with Theorems\,\ref{Sule94:Th:1} and~\ref{Sule94:Th:4}, we observe the
following:~(\romref{Th:3.3:ii}) and~(\romref{Th:3.3:iii}) can be
considered as generalizations of Theorems\,\ref{Sule94:Th:1}
and~\ref{Sule94:Th:4}, respectively.  For~(\romref{Th:3.3:ii}), we do
not assume as mentioned earlier that the prime spectrum of~${\cal A}$ is
irreducible and the plant~$P$ is strictly causal.  The rank conditions
of~${\cal T}_P$ and~${\cal W}_P$ are deleted.  For~(\romref{Th:3.3:iii}), the
commutative ring~${\cal A}$ is not restricted to a~unique factorization
domain.  The elementary factors are replaced by the generalized
elementary factors.  Although two matrices~$T$ and~$W$ in
Definition\,\ref{D:EF} are used to state Theorem\,\ref{Sule94:Th:4},
only one matrix~$T$ in Definition\,\ref{D:GEF} is used
in~(\romref{Th:3.3:iii}).

We will present the proof of the theorem
in~\S\,\ref{S:ProofMainResults}.

To make the notion of the generalized elementary factors familiar, we
present here an example of the generalized elementary factors.

\begin{example}\label{Ex:Anantharam:3.4}
On some synchronous high-speed electronic circuits such as computer
memory devices, they cannot often have nonzero small delays (for
example~\cite{bib:greenfiel85a}).  We suppose here that the system
cannot have the unit delay as a~nonzero small delay.  Further we
suppose that the impulse response of a~transfer function being stable
is finitely terminated.  Thus the set~${\cal A}$ becomes the set of
polynomials generated by~$z^2$ and~$z^3$, that is, ${\cal A}=\IR[z^2,z^3]$,
where~$z$ denotes the unit delay operator.  Then~${\cal A}$ is not a~unique
factorization domain but a~Noetherian domain.  The total field~${\cal F}$ of
fractions of~${\cal A}$ is $\IR(z^2,z^3)$, which is equal to $\IR(z)$.  The
ideal~${\cal Z}$ used to define the causality is given as the set of
polynomials in $\IR[z^2,z^3]$ whose constant terms are zero; that is,
${\cal Z}=z^2{\cal A}+z^3{\cal A}=\{ az^2+bz^3\,|\, a,b\in{\cal A}\}$.
Thus the set of causal transfer functions~${\cal P}$ is given as $n/d$,
where $n,d$ are in~${\cal A}$ and the constant term of~$d$ is nonzero; that
is,
  ${\cal P}=\{ n/(a+bz^2+cz^3)\,|\, n\in{\cal A},{\ \ } a\in\IR\backslash\{ 0\},{\ \ } b,c\in{\cal A}\}$.
Further the set of strictly causal transfer functions ${\cal P}_
{\textrm{s}}$ is given as 
  ${\cal P}_{\textrm{s}}=
     \{(a_1z^2+b_1z^3)/(a_2+b_2z^2+c_2z^3)\,|\, a_2\in\IR\backslash\{ 0\},{\ \ } a_1,b_1,b_2,c_2\in{\cal A}\}$

Since some factorized polynomials are sometimes expressed more
compactly and easier to understand than the expanded ones, we here
introduce the following notation: a~polynomial in $\IR[z]$ surrounded
by~``$\langle$'' and~``$\rangle$'' indicates that it is in~${\cal A}$ as well as in
$\IR[z]$ even though some factors between~``$\langle$'' and~``$\rangle$'' may
not be in~${\cal A}$.

Let us consider the plant below:
\begin{equation}\label{E:3.2}
  P:=\bmatrix{ (1-z^3)/(1-z^2) \cr
               (1-8z^3)/(1-4z^2) } \in{\cal P}^{2\times 1}. 
\end{equation}
The representation of the plant is not unique.  For example, the
(1,1)-entry of the plant has an alternative form $(1+z^2+z^4)/(1+z^3)$
different from the expression in (\ref{E:3.2}).  Consider
parameterizing the representation of the plant.  To do so we consider
the plant~$P$ over $\IR(z)$ rather than over~${\cal F}$.  Thus~$P$ can be
expressed as
\begin{equation}\label{E:3.3}
  P=\bmatrix{ (1+z+z^2)/(1+z) \cr
              (1+2z+4z^2)/(1+2z) } \mbox{~~~over $\IR(z)$}. 
\end{equation}
However the coefficients of all numerators and denominators in
(\ref{E:3.3}) of~$z$ with degree~$1$ are not zero.  To make them zero,
we should multiply them by $(a_1(1-z)+b_1z^2+c_1z^3)$ or
$(a_2(1-2z)+b_2z^2+c_2z^3)$ with $a_1,b_1,c_1,a_2,b_2,c_2\in{\cal A}$ as
follows
\begin{equation}\label{E:31.Oct.99.143558}
  P=\bmatrix{ 
              \frac{\langle(1+z+z^2)(a_1(1-z)+b_1z^2+c_1z^3)\rangle}%
                   {\langle(1+z)(a_1(1-z)+b_1z^2+c_1z^3)\rangle} \cr
              \frac{\langle(1+2z+4z^2)(a_2(1-2z)+b_2z^2+c_2z^3)\rangle}%
                   {\langle(1+2z)(a_2(1-2z)+b_2z^2+c_2z^3)\rangle}}.
\end{equation}
Every expression of the plant is given in the form of
(\ref{E:31.Oct.99.143558}) with $a_1,b_1,c_1,a_2,b_2,c_2$ in~${\cal A}$
provided that the denominators are not zero.  From this, we have two
observations.  One is that the plant~$P$ does not have its right- and
left-coprime factorizations over~${\cal A}$ (even so, it will be shown later
that the plant is stabilizable).  The other is that the elementary
factor of this plant cannot be consistently defined over~${\cal A}$.  Thus
we employ the notion of the generalized elementary factor.

In the following, we calculate the generalized elementary factors of
the plant.  We choose the following matrices as~$N$, $D$, and~$T$ used
in Definition\,\ref{D:GEF}:
\begin{eqnarray*}
   \bmatrix{ n_1 \cr n_2}:=
   N:=\bmatrix{(1-z^3)(1-4z^2)\cr
              (1-8z^3)(1-z^2)},\\
   \bmatrix{ d } :=
   D:=\bmatrix{(1-z^2)(1-4z^2)},{\ \ }
   T:=\bmatrix{ N^t & D^t}^t.
\end{eqnarray*}
Since $m=1$ (the number of inputs) and $n=2$ (the number of outputs),
the set~${\cal I}$ is given as
  ${\cal I}=\{\{ 1\},\{ 2\},\{ 3\}\}$
and we let
   $I_1=\{ 1\}$, 
   $I_2=\{ 2\}$,
   $I_3=\{ 3\}$. 

Let us calculate the generalized elementary factor $\Lambda_{P\!I_1}$.  Let
$i_1=1$ so that $I_1=\{ i_1\}$.  Then the $(1,i_1)$-entry of the matrix
$\Delta_{I_1}$ is~$1$ and the other entries are zero.  Thus we have $\Delta_
{I_1}=\bmatrix{1 & 0 & 0}$.  The generalized elementary factor $\Lambda_
{P\!I_1}$ is originally given as follows:
\begin{eqnarray}
  \Lambda_{P\!I_1}&=&\{\lambda\in{\cal A}\,|\,\exists K\in{\cal A}^{(m+n)\times m}\  \lambda T=K \Delta_{I_1} T\} \nonumber\\
     &=&\{\lambda\in{\cal A}\,|\,\exists k_1,k_2\in{\cal A}\ \  
             \lambda\bmatrix{ n_2 & d }^t= n_1 \bmatrix{k_1 & k_2}^t\}.
                                                \label{E:30.Oct.99.192139:1}
\end{eqnarray}
Consider (\ref{E:30.Oct.99.192139:1}) over $\IR[z]$ instead of~${\cal A}$.
Then the matrix equation in the set of (\ref{E:30.Oct.99.192139:1})
can be expressed as
\begin{eqnarray}
\lefteqn{     \lambda\bmatrix{ (1-z)(1+z)(1-2z)(1+2z+4z^2)\cr
                          (1-z)(1+z)(1-2z)(1+2z)}=}\nonumber\\
&&\hspace*{5em}   (1-z)(1-2z)(1+2z)(1+z+z^2)
              \bmatrix{k_1\cr k_2}.\label{E:31.Oct.99.144437}
\end{eqnarray}
The set of $\lambda$'s such that there exist $k_1,k_2\in\IR[z]$ satisfying
(\ref{E:31.Oct.99.144437}) is given as
  $\{(1+2z)(1+z+z^2)a\,\,|\,\,a\in\IR[z]\}$, 
denoted by~$L_1$.  Then the intersection of~$L_1$ and~${\cal A}$ is given as
follows:
\begin{equation}\label{E:09.Nov.99.103709}
   L_1\cap{\cal A}=\{\langle(1+2z)(1+z+z^2)(a(1-3z)+bz^2+cz^3)\rangle\in{\cal A}\,\,|\,\,a,b,c\in{\cal A}\}
\end{equation}
This is equal to $\Lambda_{P\!I_1}$ as shown below.  First it is obvious
that $L_1\cap{\cal A}\supset \Lambda_{P\!I_1}$.  For each
$(1+2z)(1+z+z^2)(a(1-3z)+bz^2+cz^3)$ with $a,b,c\in{\cal A}$, we have~$k_1$
and~$k_2$ as follows from (\ref{E:31.Oct.99.144437}):
\begin{eqnarray*}
  k_1&=&(1+z)(1+2z+4z^2)(a(1-3z)+bz^2+cz^3),\\
  k_2&=&(1+z)(1+2z)(a(1-3z)+bz^2+cz^3).
\end{eqnarray*}
Both~$k_1$ and~$k_2$ are in~${\cal A}$.  Hence $L_1\cap{\cal A}\subset \Lambda_{P\!I_1}$ and so
$L_1\cap{\cal A}=\Lambda_{P\!I_1}$.  By (\ref{E:09.Nov.99.103709}), we can also
consider that $\Lambda_{P\!I_1}$ is generated by $\langle(1+2z)(1+z+z^2)(1-3z)\rangle
$, $\langle(1+2z)(1+z+z^2)z^2\rangle$, and $\langle(1+2z)(1+z+z^2)z^3\rangle$.

Analogously, we can calculate the generalized elementary factors $\Lambda_
{P\!I_2}$ and $\Lambda_{P\!I_3}$ of the plant with respect to~$I_2$
and~$I_3$ as follows:
\begin{eqnarray*}
  \Lambda_{P\!I_2}
        &=&\{\langle(1+z)(1+2z+4z^2)(a(1-3z)+bz^2+cz^3)\rangle\in{\cal A}\,\,|\,\,a,b,c\in{\cal A}\},\\
  \Lambda_{P\!I_3}
        &=&\{\langle(1+z)(1+2z)(a(1-3z)+bz^2+cz^3)\rangle\in{\cal A}\,\,|\,\,a,b,c\in{\cal A}\}.
\end{eqnarray*}
Observe now that
\begin{eqnarray*}
  && \Lambda_{P\!I_1}\ni\langle(1+2z)(1-3z)(1+z+z^2)\rangle=:\lambda_{0I_1},\\
  && \Lambda_{P\!I_2}\ni\langle(1+z)(1+2z+4z^2)(1-3z+z^2)\rangle=:\lambda_{0I_2}
\end{eqnarray*}
and further
\[
    \alpha_{I_1}\lambda_{0I_1}+\alpha_{I_2}\lambda_{0I_2}=1,
\]
where
\begin{eqnarray*}
  \alpha_{I_1}&=&
\textstyle
    \frac{-4233-23646z^2-39836z^3-201780z^4-113016z^5+75344z^6}{5852}\in{\cal A},
\\
  \alpha_{I_2}&=&
\textstyle
    \frac{10085+ 18418z^2+121140z^3+131852z^4+113016z^5}{5852}\in{\cal A}.
\end{eqnarray*}
Now let
\begin{equation}\label{E:12.Nov.99.131258}
  \lambda_{I_1}:=\alpha_{I_1}\lambda_{0I_1}\in \Lambda_{P\!I_1},{\ \ }
  \lambda_{I_2}:=\alpha_{I_2}\lambda_{0I_2}\in \Lambda_{P\!I_2}.
\end{equation}
Thus $\Lambda_{P\!I_1}+\Lambda_{P\!I_2}={\cal A}$ and $\lambda_{I_1}+\lambda_{I_2}=1$.  Hence by
Theorem\,\ref{Th:3.3}, the plant~$P$ is stabilizable.
\qquad\endproof
\end{example}

\section{Intermediate Results}\label{S:IR} In this section we provide
intermediate results which will be used in the proof of our main
theorem stated in the preceding section.  This section consists of
three parts.  We first show that a~number of modules generated from
plants, controllers, and feedback systems are isomorphic to one
another.  Next we develop the results which will help to show the
existence of a~well-defined stabilizing controller.  We then give the
coprime factorizability of the plant over~${\cal A}_f$, where~$f$ is an
element of the generalized elementary factor of the plant.

\paragraph{Relationship in terms of Modules between Transfer
Matrices~{\boldmath~$P$, $C$, and $H(P,C)$}} The first intermediate
result is the relations, expressed in terms of modules, among the
matrices~$P$, $C$, and~$H(P,C)$ as well as their transposed matrices.
A~number of modules are isomorphic to one another as follows.
\begin{proposition}\label{P:4.1}
Suppose that~$P$ and~$C$ are matrices over~${\cal F}({\cal R})$.  Suppose further that
$\det(E_n+PC)$ is a~unit of ${\cal F}({\cal R})$.
\begin{romannum}
\item\label{P:4.1:i}
The following~${\cal R}$-modules are isomorphic to one another:
\begin{quote}
\begin{quote}
 ${\cal T}_{P,{\cal R}}\oplus{\cal T}_{C,{\cal R}}$,
 ${\cal T}_{H(P,C),{\cal R}}$,
 ${\cal T}_{H(P^t,C^t)^t,{\cal R}}$,
 ${\cal W}_{H(P^t,C^t),{\cal R}}$,
 ${\cal T}_{H(C,P),{\cal R}}$.
\end{quote}
\end{quote}

\item\label{P:4.1:ii}
The following~${\cal R}$-modules are isomorphic to one another:
\begin{quote}
\begin{quote}
 ${\cal W}_{P,{\cal R}}\oplus{\cal W}_{C,{\cal R}}$,
 ${\cal W}_{H(P,C),{\cal R}}$,
 ${\cal W}_{H(P^t,C^t)^t,{\cal R}}$,
 ${\cal T}_{H(P^t,C^t),{\cal R}}$,
 ${\cal W}_{H(C,P),{\cal R}}$.
\end{quote}
\end{quote}
Further for a~matrix~$X$ over ${\cal F}({\cal R})$,
\item\label{P:4.1:iii}
 ${\cal T}_{X,{\cal R}}\simeq{\cal W}_{X^t,{\cal R}}$ and  ${\cal W}_{X,{\cal R}}\simeq{\cal T}_{X^t,{\cal R}}$.
\end{romannum}
\end{proposition}

Note here that in the proposition above, the controller~$C$ need not
be a~stabilizing controller.  For the case where~$C$ is a~stabilizing
controller, see Lemma\,2 of~\cite{bib:sule94a}.

We can consider that the proposition above, especially the relations $
{\cal T}_{P,{\cal R}}\oplus{\cal T}_{C,{\cal R}}\simeq{\cal T}_{H(P,C),{\cal R}}\simeq{\cal T}_{H(C,P),{\cal R}}$, gives an
interpretation of the structure of the feedback system in the sense
that the module generated by the feedback system is given as the
direct sum of the modules generated by the plant and the controller.
In the proof (``(\romref{Th:3.3:i})$\rightarrow$(\romref{Th:3.3:ii})'') of
Theorem\,\ref{Th:3.3}, this proposition will play a~key role.
\begin{proof}
We first prove~(\romref{P:4.1:iii}).  Let~$A$ and~$B$ be matrices
over~${\cal R}$ with~$X=AB^{-1}$.  Then we have ${\cal T}_{X,{\cal R}}\simeq
M_r(\bmatrix{A^t&B^t}^t)\simeq M_c(\bmatrix{A^t&B^t})\simeq{\cal W}_{(B^{-1})^tA^t,{\cal R}}
\simeq{\cal W}_{X^t,{\cal R}}$.  The other relation ${\cal W}_{X,{\cal R}}\simeq{\cal T}_{X^t,{\cal R}}$ can be
proved in~a similar way.

Next we prove~(\romref{P:4.1:i}).  Suppose that $\det(E_n+PC)$
is a~unit of ${\cal F}({\cal R})$.  We prove the following relations in order:
  (a)~${\cal T}_{P,{\cal R}}\oplus{\cal T}_{C,{\cal R}}\simeq{\cal T}_{H(P,C),{\cal R}}$, 
  (b)~${\cal T}_{H(P,C),{\cal R}}\simeq{\cal T}_{H(P^t,C^t)^t,{\cal R}}$, 
  (c)~${\cal T}_{H(P^t,C^t)^t,{\cal R}}\simeq{\cal W}_{H(P^t,C^t),{\cal R}}$, 
  (d)~${\cal T}_{H(P,C),{\cal R}}\simeq{\cal T}_{H(C,P),{\cal R}}$. 

(a) of~(\romref{P:4.1:i}).  The proof of (a) follows mainly the proof
of Lemma\,2 in~\cite{bib:sule94a}.  By virtue of Lemma\,\ref{L:2.1},
it is sufficient to show the relation ${\cal T}_{P,{\cal R}}\oplus{\cal T}_{C,{\cal R}}\simeq
M_r(\bmatrix{N_H^t&d_HE_{m+n}}^t)$ with $N_H\in({\cal R})_{m+n}$, $d_H\in{\cal R}$,
where $H(P,C)=N_Hd_H^{-1}$.
Let~$N$, $N_C$ be matrices over~${\cal R}$ and~$d$, $d_C$ be in~${\cal R}$ with
$P=Nd^{-1}$ and $C=N_Cd_C^{-1}$.  Further, let
\[
  Q=\bmatrix{  d_CE_n & N \cr
                -N_C    & dE_m   },{\ \ }
  S=\bmatrix{  d_CE_n & O \cr
                  O      & dE_m   }.
\]
From these we have
 ${\cal T}_{P,{\cal R}}\oplus{\cal T}_{C,{\cal R}}\simeq M_r(\bmatrix{ Q^t & S^t}^t)$.
In addition, since $\det(E_n+PC)$ is a~unit of ${\cal F}({\cal R})$, the
matrix~$N_H$ is nonsingular,
so that 
  $M_r(\bmatrix{ Q^t & S^t}^t)\simeq M_r(\bmatrix{ Q^t & S^t}^t (\det(N_H)E_{m+n}))$ 
holds.  A~simple calculation shows that
\[
     \bmatrix{ Q \cr
               S} (\det(N_H)E_{m+n})
   = \bmatrix{ d_HE_{m+n}\cr
               N_H }\adj(N_H)S.
\]
Because both matrices~$S$ and $\adj(N_H)$ are nonsingular, we finally have
that
\[
{\cal T}_{P,{\cal R}}\oplus{\cal T}_{C,{\cal R}} \simeq M_r(\bmatrix{  Q \cr
                             S})
          \simeq M_r(\bmatrix{  Q \cr
                           S} (\det(N_H)E_{m+n}))
          \simeq M_r(\bmatrix{ d_HE_{m+n} \cr
                           N_H })
          \simeq {\cal T}_{H(P,C),{\cal R}}.
\]

(b) of~(\romref{P:4.1:i}). 
Observe that the following relation holds:
\begin{equation}\label{E:13.Nov.97.181751}
  H(P^t,C^t)^t=\bmatrix{ O   & E_m \cr
                         E_n & O}
               H(P,C)
               \bmatrix{ O   & E_n \cr
                         E_m & O}.
\end{equation}
Let~$N_H$ and~$d_H$ be a~matrix over~${\cal R}$ and a~scalar of~${\cal R}$,
respectively, with $H(P,C)=N_Hd_H^{-1}$.  Then
(\ref{E:13.Nov.97.181751}) can be rewritten as
follows:
\[%
  H(P^t,C^t)^t=\bmatrix{ O   & E_m \cr
                         E_n & O}
                N_H
               (\bmatrix{ O   & E_m \cr
                          E_n & O}d_H)^{-1}.
\]%
Hence, we have matrices~$A$ and~$B$ over~${\cal R}$ with $H(P^t,C^t)^t=AB^{-1}$
such that
\[%
  A=\bmatrix{ O & E_m \cr E_n & O} N_H,{\ \ }
  B=\bmatrix{ O & E_m \cr E_n & O}d_H.
\]%
This gives the relation ${\cal T}_{H(P,C),{\cal R}}\simeq{\cal T}_{H(P^t,C^t)^t,{\cal R}}$.

(c)~of~(\romref{P:4.1:i}).  This is directly
obtained by applying~(iii) to the matrix $H(P^t,C^t)^t$.

(d)~of~(\romref{P:4.1:i}).  Between the matrices $H(P,C)$ and
$H(C,P)$, the following relation holds:
\[%
     H(C,P)=\bmatrix{O&-E_n\cr
                     E_m&O}H(P,C)
            \bmatrix{O&E_m\cr
                     -E_n&O}.
\]%
Letting~$N_H$ and~$d_H$ be a~matrix over~${\cal R}$ and a~scalar of~${\cal R}$
with $H(P,C)=N_Hd_H^{-1}$ as in (b), we have matrices~$N_H'$
and~$D_H'$ over~${\cal R}$ such that
\[
  \bmatrix{N_H'\cr D_H'}
 = \left[
     \begin{array}{cc}
       \begin{array}{cc}
            O&-E_n\\
            E_m&O
       \end{array} & 
       \mbox{\Large$O$}\\
       \mbox{\Large$O$} &
       \begin{array}{cc}
            O&-E_n\\
            E_m&O
       \end{array}
     \end{array}
   \right]
  \bmatrix{N_H\cr d_H E_{m+n}}
\]
holds.  Since $H(C,P)=N_H'{D_H'}^{-1}$ holds and the first matrix of
the right-hand side of the equation above is invertible, the relation
${\cal T}_{H(P,C),{\cal R}}\simeq{\cal T}_{H(C,P),{\cal R}}$ holds.

Finally, arguments similar to~(\romref{P:4.1:i})
prove~(\romref{P:4.1:ii}). 
\qquad
\end{proof}

Before moving to the next, we prove an easy lemma useful to give
results for the transposed plants.
\begin{lemma}\label{L:4.2}
A~plant~$P$ is stabilizable if and only if its transposed plant~$P^t$
is.  Moreover, in the case where the plant~$P$ is stabilizable, $C$ is
a~stabilizing controller of~$P$ if and only if~$C^t$ is a~stabilizing
controller of the transposed plant~$P^t$.
\end{lemma}
\begin{proof}
(Only If)~Suppose that a~plant~$P$ is stabilizable.  Let~$C$ be
a~stabilizing controller of~$P$.  First, $(P^t,C^t)$ is in $\Fad$,
since $(P,C)\in\Fad$ and $\det(E_n+PC)=\det(E_m+P^tC^t)$.
From~(\ref{E:13.Nov.97.181751}) in the proof of
Proposition\,\ref{P:4.1}, if $H(P,C)\in({\cal A})_{m+n}$, then $H(P^t,C^t)\in
({\cal A})_{m+n}$.

(If) Because $(P^t)^t=P$, the ``If'' part can be proved analogously.
\qquad
\end{proof}

\paragraph{${\cal Z}$-nonsingularity of Transfer Matrices}
In order to prove the stabilizability of the given causal plant, which
will be necessary in the proof of the main theorem
(Theorem\,\ref{Th:3.3}), we should show the existence of the
stabilizing controller.  To do so, we will need to show that the
denominator matrix of the stabilizing controller is nonsingular.  The
following result will help this matter.

\begin{lemma}\label{L:4.3}
Suppose that there exist matrices~$A$, $B$, $C_1$, $C_2$ over~${\cal A}$
such that the following square matrix is ${\cal Z}$-nonsingular:
\begin{equation}\label{E:L:4.3:1}
   \bmatrix{ A&C_1 \cr
             B&C_2},
\end{equation}
where the matrix~$A$ is square and the matrices~$A$ and~$B$ have same
number of columns.  Then there exists a~matrix~$R$ over~${\cal A}$ such that
the matrix $A+RB$ is ${\cal Z}$-nonsingular.
\end{lemma}

Before starting the proof, it is worth reviewing some easy facts about
the prime ideal~${\cal Z}$.
\begin{remark}\label{R:4.4}
(i) If~$a$ is in~${\cal A}\backslash{\cal Z}$ and expressed as $a=b+c$ with
$b,c\in{\cal A}$, then at least one of~$b$ and~$c$ is in~${\cal A}\backslash{\cal Z}$.
(ii)
If~$a$ is in ${\cal A}\backslash{\cal Z}$ and~$b$ is in~${\cal Z}$, then the sum $a+b$ is in
${\cal A}\backslash{\cal Z}$.
(iii) Every factor in~${\cal A}$ of an element of ${\cal A}\backslash{\cal Z}$ belongs to ${\cal A}\backslash
{\cal Z}$ (that is, if $a,b\in{\cal A}$ and $ab\in{\cal A}\backslash{\cal Z}$, then $a,b\in{\cal A}\backslash{\cal Z}$).
\qquad\endproof
\end{remark}

They will be used in the proofs of Lemma\,\ref{L:4.3} and
Theorem\,\ref{Th:3.3}.
\begin{proofof}{Lemma\,\ref{L:4.3}}
This proof mainly follows that of~Lemma\,4.4.21
of~\cite{bib:vidyasagar85a}.

If the matrix~$A$ itself is ${\cal Z}$-nonsingular, then we can select the
zero matrix as~$R$.  Hence we assume in the following that~$A$ is ${\cal Z}
$-singular.

Since~(\ref{E:L:4.3:1}) is ${\cal Z}$-nonsingular, there exists a~full-size
minor of $\bmatrix{A^t & B^t}^t$ in ${\cal A}\backslash{\cal Z}$ by Laplace's expansion
of~(\ref{E:L:4.3:1}) and Remark\,\ref{R:4.4}(i,iii).  Let~$a$ be such
a ${\cal Z}$-nonsingular full-size minor of
  $\bmatrix{A^t & B^t}^t$
having as few rows from~$B$ as possible.  

We here construct a~matrix~$R$ such that
  $\det(A+RB)=\pm a+z$
with~$z\in{\cal Z}$.  Since~$A$ is ${\cal Z}$-singular, the full-size minor~$a$
must contain at least one row of~$B$ from the matrix $\bmatrix{A^t &
B^t}^t$.  Suppose that~$a$ is obtained by excluding the rows
$i_1,\ldots, i_k$ of~$A$ and including the rows $j_1,\ldots, j_k$ of~$B$,
where both of $i_1,\ldots, i_k$ and $j_1,\ldots, j_k$ are in ascending order.
{Now} define $R=(r_{ij})$ by
  $r_{i_1 j_1}=\cdots=r_{i_k j_k}=1$
and $r_{ij}=0$ for all other~$i$,~$j$.  Observe that
$\det(A+RB)$ is expanded in terms of full-size
minors of the matrices $\bmatrix{E & R}$ and
  $\bmatrix{A^t & B^t}^t$
from the factorization
  $A+RB=\bmatrix{E & R}
                       \bmatrix{A^t & B^t}^t$ 
by the Binet-Cauchy formula.  Every minor of $\bmatrix{E & R}$
containing more than~$k$ columns of~$R$ is zero.  By the method of
choosing the rows from 
  $\bmatrix{A^t & B^t}^t$
for the full-size minor~$a$, every full-size minor of
  $\bmatrix{A^t & B^t}^t$
having less than~$k$ rows of~$B$ is in~${\cal Z}$.  There is only one
nonzero minor of $\bmatrix{E & R}$ containing exactly~$k$ columns
of~$R$, which is obtained by excluding the columns $i_1,\ldots, i_k$ of
the identity matrix~$E$ and including the columns $j_1,\ldots, j_k$
of~$R$; it is equal to~$\pm 1$.  From the Binet-Cauchy formula, the
corresponding minor of $\bmatrix{A^t & B^t}^t$ is~$a$.  As a~result,
$\det(A+RB)$ is given as a {sum} of~$\pm a$ and elements in~${\cal Z}$.  By
Remark\,\ref{R:4.4}(ii), the sum is in ${\cal A}\backslash{\cal Z}$ and so is
$\det(A+RB)$.  The matrix $A+RB$ is now ${\cal Z}$-nonsingular.  \qquad
\end{proofof}

\paragraph{Coprimeness and Generalized Elementary Factors}
We present here that for each non-nilpotent element~$\lambda$ of the
generalized elementary factors, the plant has a~right-coprime
factorization over~${\cal A}_{\lambda}$.  This will be independent of the
stabilizability of the plant.

\begin{lemma}\label{L:4.5}
\textup{(cf.\,Proposition\,2.2 of~\cite{bib:mori97b})~} Let $\Lambda_
{P\!I}$ be the generalized elementary factor of the plant~$P$ with
respect to~$I\in{\cal I}$ and further let $\sqrt{\Lambda_{P\!I}}$ denote the
radical of $\Lambda_{P\!I}$ (as an ideal).  Suppose that~$\lambda$ is in $\sqrt
{\Lambda_{P\!I}}$ but not nilpotent.  Then, the ${\cal A}_{\lambda}$-module ${\cal T}_{P,{\cal A}_{\lambda}%
}$ is free of rank~$m$.
\end{lemma}
\begin{proof}
Let~$T,N,D$ be the matrices over~${\cal A}$ as in Definition\,\ref{D:GEF}.
Recall that~${\cal T}_{P,{\cal A}_{\lambda}}$ denotes the~${\cal A}_{\lambda}$-module generated by
rows of the matrix~$T$.  By Definition\,\ref{D:GEF}, there exists
a~matrix~$K$ over~${\cal A}$ such that $\lambda^rT=K \Delta_I T$ holds for some
positive integer~$r$.  Then we have a~factorization of the matrix~$T$
over~${\cal A}_{\lambda}$ as
  $T=(\lambda^{-r}K)(\Delta_I T)$,
where all entries of the matrix $\lambda^{-r}K$ are in~${\cal A}_{\lambda}$.
In order to show that the~${\cal A}_{\lambda}$-module ${\cal T}_{P,{\cal A}_{\lambda}}$ is free of
rank~$m$, provided that~$\lambda$ is not nilpotent, it is sufficient to
prove the following two facts: (i) The matrix $\Delta_I T$ is nonsingular
over~${\cal A}_{\lambda}$.~~(ii) There is a~matrix~$X$ such that the matrix
$\bmatrix{\lambda^{-r}K&X}$ is invertible over~${\cal A}_{\lambda}$ and the matrix
equation
  $T=\bmatrix{\lambda^{-r}K&X} \bmatrix{ (\Delta_I T)^t & O}^t$ 
holds.

(i).~Observe that the matrix~$D$ is nonsingular over~${\cal A}_{\lambda}$ as well
as over~${\cal A}$.  Since $D=\Delta_{\{ n+1,\ldots, m+n\}}T=(\lambda^{-r}\Delta_
{\{ n+1,\ldots, m+n\}}K)(\Delta_I T)$ holds (note that the suffix of the
symbol~$\Delta$ is an ordered set of~$m$ distinct integers between~$1$ and
$m+n$ as before Definition\,\ref{D:GEF}), the matrix $\Delta_I T$ is also
nonsingular over~${\cal A}_{\lambda}$ provided that~$\lambda$ is not nilpotent.

(ii).~Let $\overline{i_1},\ldots,\overline{i_n}$ be~$n$ distinct integers
in ascending order between~$1$ and $m+n$ excluding the integers
in~$I$.  Then let~$X$ be the matrix whose $(\overline{i_k},k)$-entry
is~$1$ for each~$\overline{i_k}$ and zero otherwise.
Then the determinant of $\bmatrix{\lambda^{-r}K&X}$ becomes~$\pm 1$ since the
matrix $\lambda^{-r}\Delta_IK$ is the identity matrix of $({\cal A}_{\lambda})_m$.
\qquad
\end{proof}

The lemma below will be used in the proof (``(\romref{Th:3.3:ii})$\rightarrow
$(\romref{Th:3.3:iii})'') of the main theorem by letting ${\cal R}={\cal A}_f$,
where~$f$ is an element of the generalized elementary factor of the
plant but not nilpotent.

\begin{lemma}\label{L:4.6}
If~${\cal R}$-module~${\cal T}_{P,{\cal R}}$ \,$\Bigl({\cal W}_{P,{\cal R}}\Bigr)$ is free of
rank~$m$ \,$\Bigl(n\Bigr)$, there exist matrices~$A$ and~$B$
$\Bigl(\widetilde{A}$ and $\widetilde{B}\Bigr)$ over~${\cal R}$ such that
$(A,B)$ is a~right-coprime factorization
$\Bigl((\widetilde{A},\widetilde{B})$ is a~left-coprime
factorization$\Bigr)$ of the plant~$P$ ($\in{\cal F}({\cal R})^{n\times m}$) over~${\cal R}$.
\end{lemma}
\begin{proof}
This lemma is an analogy of the result given in the proof of
Lemma\,3 of~\cite{bib:sule94a}.
See the proof of Lemma\,3
of~\cite{bib:sule94a}.
\qquad
\end{proof}

\begin{example}\label{Ex:12.Nov.99.111204}
(Continued) Here we continue Example\,\ref{Ex:Anantharam:3.4}.  Let us
follow Lemmas\,\ref{L:4.5} and~\ref{L:4.6} with the plant of
(\ref{E:3.2}).  Let the notation be as in
Example\,\ref{Ex:Anantharam:3.4}.

First we proceed along the proof of Lemma\,\ref{L:4.5}.  As an
example, we pick $I_1\in{\cal I}$ as~$I$ and $\lambda_{I_1}\in \Lambda_{P\!I_1}$ as~$\lambda$
in the proof of Lemma\,\ref{L:4.5}.  Recall that for each $\lambda\in \Lambda_
{P\!I}$, there exists a~matrix~$K$ such that $\lambda T=K \Delta_IT$ holds.  In
the case of $\lambda_{I_1}\in \Lambda_{P\!I_1}$, the matrix~$K$ is given as
\begin{equation}\label{E:11.Nov.99.213238}
  K=\bmatrix{k_1\cr k_2 \cr k_3}
   =\bmatrix{ \lambda_{I_1} \cr
              \alpha_{I_1}\langle(1+z)(1-3z)(1+2z+4z^2)\rangle \cr
              \alpha_{I_1}\langle(1+z)(1+2z)(1-3z)\rangle }.
\end{equation}
Thus we have
the factorization $T=(\lambda^{-r}K)(\Delta_{I_1} T)$:
\[
   \bmatrix{(1-z^3)(1-4z^2)\cr
            (1-8z^3)(1-z^2)\cr
            (1-z^2)(1-4z^2)}
  =
   \bmatrix{ 1  \cr
             \lambda_{I_1}^{-1}\alpha_{I_1}\langle(1+z)(1-3z)(1+2z+4z^2)\rangle \cr
             \lambda_{I_1}^{-1}\alpha_{I_1}\langle(1+z)(1+2z)(1-3z)\rangle}
   \bmatrix{(1-z^3)(1-4z^2)},
\]
where $r=1$ and $\Delta_{I_1}T=\bmatrix{(1-z^3)(1-4z^2)}$.  As shown in
part~(i) of the proof of Lemma\,\ref{L:4.5}, $\Delta_
{I_1}T=\bmatrix{(1-z^3)(1-4z^2)}$ is nonsingular.

The matrix~$X$ in part~(ii) of the proof of Lemma\,\ref{L:4.5} is
given as $X=[{ 0~1~0 \atop 0~0~1}]^t$ by letting $\overline{i}_1=2$
and $\overline{i}_2=3$ according to $I_1=\{ 1\}$.  We can see that the
matrix
\begin{equation}\label{E:02.Nov.99.132204}
  \bmatrix{\lambda^{-1}K & X}
   =\bmatrix{ 1                                         & 0 & 0\cr
              \lambda_{I_1}^{-1}\alpha_{I_1}\langle(1+z)(1-3z)(1+2z+4z^2)\rangle & 1 & 0\cr
              \lambda_{I_1}^{-1}\alpha_{I_1}\langle(1+z)(1+2z)(1-3z)\rangle        & 0 & 1}
\end{equation}
is invertible.  Therefore the ${\cal A}_{\lambda_{I_1}}$-module ${\cal T}_{P,{\cal A}_{\lambda_
{I_1}}}$ is free and its rank is~$1$.  (However we will see that the $
{\cal A}$-module~${\cal T}_P$ is not free, see Example\,\ref{Ex:Anantharam:5.3})

From (\ref{E:02.Nov.99.132204}) and the matrix equation $T=\bmatrix{\lambda
^{-r}K&X} \bmatrix{ (\Delta_{I_1} T)^t & O}^t$, we let
\begin{eqnarray}
  \bmatrix{N_{I_1} \cr D_{I_1}} &:=&
    \bmatrix{ 1  \cr
              \lambda_{I_1}^{-1}\alpha_{I_1}\langle(1+z)(1-3z)(1+2z+4z^2)\rangle \cr
              \lambda_{I_1}^{-1}\alpha_{I_1}\langle(1+z)(1+2z)(1-3z)\rangle}
    =\lambda_{I_1}^{-1}K,\label{E:02.Nov.99.150110:1}\\
 \bmatrix{\widetilde{Y}_{I_1} & \widetilde{X}_{I_1}\cr
          \times & \times}
   &:=&\bmatrix{ 1                                         & 0 & 0\cr
                -\lambda_{I_1}^{-1}\alpha_{I_1}\langle(1+z)(1-3z)(1+2z+4z^2)\rangle & 1 & 0\cr
                -\lambda_{I_1}^{-1}\alpha_{I_1}\langle(1+z)(1+2z)(1-3z)\rangle        & 0 & 1}
\nonumber\\
   &&~~~     =\bmatrix{\lambda_{I_1}^{-1}K & X}^{-1}, \label{E:02.Nov.99.150110:2}
\end{eqnarray}
where $N_{I_1}\in{\cal A}_{\lambda_{I_1}}^{2\times 1}$, $\widetilde{Y}_{I_1}\in{\cal A}_{\lambda_
{I_1}}^{1\times 2}$, $D_{I_1},\widetilde{X}_{I_1}\in({\cal A}_{\lambda_{I_1}})_1$, and
$\times$ denotes some matrix.  Then $(N_{I_1},D_{I_1})$ is a~right-coprime
factorization of the plant over~${\cal A}_{\lambda_{I_1}}$ with
  $\widetilde{Y}_{I_1} N_{I_1}+\widetilde{X}_{I_1} D_{I_1}=E_1$.
which is consistent with Lemma\,\ref{L:4.6}.
\qquad\endproof
\end{example}

\section{Proof of Main Results}\label{S:ProofMainResults}
{Now} we give the proof of the main theorem.
\begin{proofof}{Theorem\,\ref{Th:3.3}}
We prove the following relations in order:
  (a) ``(\romref{Th:3.3:i})$\rightarrow$(\romref{Th:3.3:ii}),''
  (b) ``(\romref{Th:3.3:ii})$\rightarrow$(\romref{Th:3.3:iii}),'' and
  (c) ``(\romref{Th:3.3:iii})$\rightarrow$(\romref{Th:3.3:i}).''

\noindent
(a) ``(\romref{Th:3.3:i})$\rightarrow$(\romref{Th:3.3:ii})'': Suppose that~$C$
is a~stabilizing controller of the plant~$P$.  Then, the ${\cal A}$-module $
{\cal T}_{H(P,C)}$ is obviously free.  By the relation ${\cal T}_{P,{\cal R}}\oplus{\cal T}_{C,{\cal R}}
\simeq{\cal T}_{H(P,C),{\cal R}}$ in Proposition\,\ref{P:4.1}(\romref{Th:3.3:i}), we
have that the ${\cal A}$-module~${\cal T}_P$ is projective.  By using
Proposition\,\ref{P:4.1}(iii) and Lemma\,\ref{L:4.2}, the projectivity
of the ${\cal A}$-module~${\cal W}_P$ can be proved analogously.

\noindent
(b) ``(\romref{Th:3.3:ii})$\rightarrow$(\romref{Th:3.3:iii})'': Suppose
that~(\romref{Th:3.3:ii}) holds, that is, the modules~${\cal T}_P$ and ${\cal W}_
P$ are projective.  We let~$T,N,D$ be the matrices over~${\cal A}$ as in
Definition\,\ref{D:GEF}.  According to Theorem\,IV.32
of~\cite[p.295]{bib:brmacdonald84a}, 
there exist finite sets~$F_1$ and~$F_2$ such that~(1) each of them
generates~${\cal A}$ and~(2) for any $f\in F_1$ $\Bigl(f\in F_2\Bigr)$, the~${\cal A}_
f$-module ${\cal T}_{P,{\cal A}_f}$ $\Bigl({\cal W}_{P,{\cal A}_ f}\Bigr)$ is free.  Let~$F$
be the set of all $f_1f_2$'s with $f_1\in F_1$ and $f_2\in F_2$.  Then~$F$
generates~${\cal A}$ again, and the~${\cal A}_f$-modules ${\cal T}_{P,{\cal A}_f}$ and ${\cal W}_{P,
{\cal A}_ f}$ are free for every $f\in F$.  We suppose without loss of
generality that the sets~$F_1$, $F_2$, and~$F$ do not contain any
nilpotent element because $1+x$ is a~unit of~${\cal A}$ for any
nilpotent~$x$ (cf.\,\cite[p.10]{bib:atiyah69a}).  (However, we note
that other zerodivisors cannot be excluded from the set~$F$.)  The
rank of the free~${\cal A}_f$-module ${\cal T}_{P,{\cal A}_f}$ is~$m$, since~$m$ rows of
the denominator matrix~$D$ are independent over~${\cal A}_f$ as well as
over~${\cal A}$.  Analogously the rank of ${\cal W}_{P,{\cal A}_f}$ is~$n$.

In order to show that~(\romref{Th:3.3:iii}) holds, it suffices to show
that the relation $\sum_{f\in F}(f^{\xi})\subset\sum_{I\in{\cal I}}\Lambda_{P\!I}$ holds for
a~sufficiently large integer~$\xi$.  Once this relation is obtained,
since $\sum_{f\in F}(f^{\xi})={\cal A}$ holds, we have $ \sum_{I\in{\cal I}}\Lambda_{P\!I}=
{\cal A}$.

Let~$f$ be an arbitrary but fixed element of~$F$.
Let~$V_f$ be a~square matrix of size~$m$ whose rows are~$m$ distinct
generators of the~${\cal A}_ f$-module $M_r(\bmatrix{ N^t&D^t}^t)$ ($\simeq{\cal T}_
{P,{\cal A}_f}$).  We assume without loss of generality that~$V_f$ is over~$
{\cal A}$, that is, the denominators of all entries of~$V_f$ are~$1$.
Otherwise if~$V_f$ is over~${\cal A}_f$ but not over~${\cal A}$, $V_f$ multiplied
by~$f^x$, with a~sufficiently large integer~$x$, will be over~${\cal A}$, so
that we can consider such $V_ff^x$ as ``$V_f$.''  
Thus the following matrix equation holds over~${\cal A}
$: 
\begin{equation}\label{E:11.Nov.99.213521}
  f^{\nu}T=\Kf V_f
\end{equation}
with a~nonnegative integer~$\nu$ and a~matrix $\Kf \in{\cal A}^{(m+n)\times m}$.

In order to prove the relation $\sum_{f\in F}(f^{\xi})\subset\sum_{I\in{\cal I}}\Lambda_
{P\!I}$, we will first show the relation
\begin{equation}\label{E:14.Feb.97.113224}
  I_{m{\cal A}}(f^{\nu}\Kf)\subset\sum_{I\in{\cal I}}\Lambda_{P\!I}
\end{equation}
and then
\begin{equation}\label{E:08.Feb.97.134530}
  (f^{\xi})\subset I_{m{\cal A}}(f^{\nu}\Kf).
\end{equation}

Observe first that $\det(f^{\nu}\Delta_I\Kf )\in \Lambda_{P\!I}$ because
\[
  \det(f^{\nu}\Delta_I\Kf )T=f^{m \nu}\Kf \adj(\Delta_I\Kf )\Delta_IT.
\]
Since every element of $I_{m{\cal A}}(f^{\nu}\Kf)$ is an~${\cal A}$-linear
combination of $\det(f^{\nu}\Delta_I\Kf)$'s for all $I\in{\cal I}$, we
have~(\ref{E:14.Feb.97.113224}).

We next present~(\ref{E:08.Feb.97.134530}).  Let~$N_0$ and~$D_0$ be
matrices with $\Kf=\bmatrix{N_0^t & D_0^t}^t$.  Since each row
of~$V_f$ is generated by rows of $\bmatrix{N^t&D^t}^t$ over~${\cal A}_f$,
there exist matrices~$\widetilde{Y}_0$ and~$\widetilde{X}_0$ over~${\cal A}_
f$ such that $V_f=\bmatrix{\widetilde{Y}_0 f^{\nu}&\widetilde{X}_0f^{\nu}}
\bmatrix{N^t&D^t}^t$.  Thus, since~$V_f$ is nonsingular over~${\cal A}_f$,
we have 
  $\bmatrix{\widetilde{Y}_0&\widetilde{X}_0}
   \bmatrix{N_0^t&D_0^t}^t=E_m$, 
which implies that $(N_0,D_0)$ is a~right-coprime factorization of the
plant~$P$ over~${\cal A}_f$.  Recall here that ${\cal W}_{P,{\cal A}_f}$ is free of
rank~$n$.  Thus by Lemma\,\ref{L:4.6} there exist matrices
$\widetilde{N}_0$ and $\widetilde{D}_0$ such that
$(\widetilde{N}_0,\widetilde{D}_0)$ is a~left-coprime factorization of
the plant~$P$ over~${\cal A}_f$.  Let~$Y_0$ and~$X_0$ be matrices over~${\cal A}_
f$ such that $\widetilde{N}_0Y_0+\widetilde{D}_0X_0=E_n$ holds.  Then
we have the following matrix equation:
\begin{equation}\label{E:11.Feb.99.105503}
   \bmatrix{  \widetilde{Y}_0 & \widetilde{X}_0 \cr
             -\widetilde{D}_0 & \widetilde{N}_0 }
   \bmatrix{  N_0 & -X_0 \cr
              D_0 &  Y_0 }
  =\bmatrix{ E_m & -\widetilde{Y}_0X_0+\widetilde{X}_0Y_0 \cr
             O   & E_n}.
\end{equation}
Denote by~$R$ the matrix $\bmatrix{-X_0^t & Y_0^t}^t$.  Then the
matrix $\bmatrix{ \Kf & R }$ is invertible over~${\cal A}_f$ since the
right-hand side of~(\ref{E:11.Feb.99.105503}) is invertible.  For each
$I\in{\cal I}$, let $\overline{I}$ be the ordered set of~$n$ distinct
integers between~$1$ and $m+n$ excluding~$m$ integers in~$I$ and let
$\overline{i_1},\ldots,\overline{i_n}$ be elements of $\overline{I}$ in
ascending order.  Let $\Delta_{\overline{I}}\in{\cal A}^{m\times(m+n)}$ denote the
matrix whose $(k,\overline{i_k})$-entry is~$1$ for $\overline{i_k}\in
\overline{I}$ and zero otherwise.  Then, by using Laplace's expansion,
the following holds:
\[
  \det(\bmatrix{\Kf & R})=\sum_{I\in{\cal I}} (\pm\det(\Delta_I\Kf) \det(\Delta_{\overline{I}}R)),
\]
which is a~unit of~${\cal A}_f$.  From this and since the ideal $I_{m{\cal A}_
f}(\Kf)$ is generated by $\det(\Delta_I\Kf)$'s for all~$I\in{\cal I}$, we have
$I_{m{\cal A}_f}(\Kf)={\cal A}_f$.  This equality over~${\cal A}_f$ gives
(\ref{E:08.Feb.97.134530}) for a~sufficiently large integer~$\xi$.

From (\ref{E:14.Feb.97.113224}) and (\ref{E:08.Feb.97.134530}), the
relation $\sum_{f\in F}(f^{\xi})\subset\sum_{I\in{\cal I}}\Lambda_{P\!I}$ holds.  Therefore
we conclude that the relation $\sum_{I\in{\cal I}}\Lambda_{P\!I}={\cal A}$ holds.

\bigskip\par\noindent (c) ``(\romref{Th:3.3:iii})$\rightarrow
$(\romref{Th:3.3:i})'': To prove the stabilizability, we will
construct a~stabilizing controller of the causal plant~$P$ from
right-coprime factorizations over~${\cal A}_f$ with some $f$'s in~${\cal A}$.
Let~$N$ and~$D$ be matrices over~${\cal A}$ such that $P=ND^{-1}$ and~$D$ is
${\cal Z}$-nonsingular.  From (\ref{E:Th:3.2:1}), there exist~$\lambda_I$'s such
that $\sum \lambda_I=1$, where~$\lambda_I$ is an element of generalized
elementary factor $\Lambda_{P\!I}$ of the plant~$P$ with respect to~$I$
in~${\cal I}$; that is, $\lambda_I\in \Lambda_{P\!I}$.  {Now} let these~$\lambda_I$'s be
fixed.  Further, let~${\cal I}^\sharp$ be the set of~$I$'s of these nonzero~$\lambda_
I$'s; that is,
  $\sum_{I\in{\cal I}^\sharp} \lambda_I=1$. 
As in (b), we can consider without loss of generality that for every
$I\in{\cal I}^\sharp$, $\lambda_I$ is not a~nilpotent element of~${\cal A}$.  For each $I\in
{\cal I}^\sharp$, the ${\cal A}_{\lambda_I}$-module ${\cal T}_{P,{\cal A}_{\lambda_I}}$ is free of rank~$m$
by Lemma\,\ref{L:4.5}.  As in~(b) again, let~$V_{\lambda_I}$ be a~square
matrix of size~$m$ whose rows are~$m$ distinct generators of the~${\cal A}_
{\lambda_I}$-module $M_r(\bmatrix{N^t&D^t}^t)$ ($\simeq{\cal T}_{P,{\cal A}_{\lambda_ I}}$) and
assume without loss of generality that~$V_{\lambda_I}$ is over~${\cal A}$.  Then
there exist matrices $\widetilde{X}_I$, $\widetilde{Y}_I$, $N_I$,
$D_I$ over~${\cal A}_{\lambda_I}$ such that
\begin{equation}\label{E:09.Feb.97.151238}
  \begin{array}{c}
    \bmatrix{N^t & D^t}^t=\bmatrix{N_I^t & D_I^t}^tV_{\lambda_I},{\ \ }
    \bmatrix{\widetilde{Y}_I & \widetilde{X}_I}\bmatrix{N^t & D^t}^t=V_{\lambda_I},\\
     \widetilde{Y}_IN_I
    +\widetilde{X}_ID_I
    =E_m
  \end{array}
\end{equation}
with $P=N_ID_I^{-1}$ over ${\cal F}({\cal A}_{\lambda_I})$.

We here present a~formula of a~stabilizing controller which is
constructed from the matrices~$\widetilde{X}_I$ and~$\widetilde{Y}_I$.
For any positive integer~$\omega$, there are coefficients~$a_I$'s in~${\cal A}$
with $ \sum_{I\in{\cal I}^\sharp} a_I \lambda_I^{\omega}=1$.  Let~$\omega$ be a~sufficiently
large integer.  Thus the matrices $\lambda_I^{\omega}D_I\widetilde{X}_I$ and $\lambda_
I^{\omega} D_I \widetilde{Y}_I$ are over~${\cal A}$ for every $I\in{\cal I}^\sharp$.  The
stabilizing controller we will construct has the form
\begin{equation}\label{E:08.Feb.97.215931}
  C=(\sum_{I\in{\cal I}^\sharp} a_I \lambda_I^{\omega} D_I \widetilde{X}_I)^{-1}
    (\sum_{I\in{\cal I}^\sharp} a_I \lambda_I^{\omega} D_I \widetilde{Y}_I). {\ \ } {\ \ }
\end{equation}

In the following we first consider that the matrix $(\sum_{I\in{\cal I}^\sharp}
a_I \lambda_I^{\omega} D_I \widetilde{X}_I)$ is ${\cal Z}$-nonsingular and show that
the plant is stabilized by the matrix~$C$ of
(\ref{E:08.Feb.97.215931}).  After showing it, we will be concerned
with the case where the matrix $(\sum_{I\in{\cal I}^\sharp} a_I \lambda_I^{\omega} D_I
\widetilde{X}_I)$ is ${\cal Z}$-singular.

Suppose that the matrix $(\sum_{I\in{\cal I}^\sharp} a_I \lambda_I^{\omega} D_I
\widetilde{X}_I)$ is ${\cal Z}$-nonsingular.  To prove that~$C$ is a~stabilizing controller of~$P$, it is sufficient to show that $(P,C)\in
\Fad$ and that four blocks of (\ref{E:H(P,C)}) are over~${\cal A}$.

We first show that $(P,C)\in\Fad$.  The following matrix equation
holds:
\begin{eqnarray*}
     E_m+CP
 &=& E_m+(\sum_{I\in{\cal I}^\sharp} a_I \lambda_I^{\omega} D_I \widetilde{X}_I)^{-1} 
         (\sum_{I\in{\cal I}^\sharp} a_I \lambda_I^{\omega} D_I \widetilde{Y}_I)
         ND^{-1}\\
 &=& (\sum_{I\in{\cal I}^\sharp} a_I \lambda_I^{\omega} D_I \widetilde{X}_I)^{-1}\\
 &&  {\ \ } {\ \ \ } {\ \ }
     \bigl(
        (\sum_{I\in{\cal I}^\sharp} a_I \lambda_I^{\omega} D_I \widetilde{X}_I)D
       +(\sum_{I\in{\cal I}^\sharp} a_I \lambda_I^{\omega} D_I \widetilde{Y}_I)N
     \bigr)
      D^{-1}.
\end{eqnarray*}
By the (1,1)-block of (\ref{E:08.Feb.97.220842}), we have
\begin{equation}\label{E:08.Feb.97.224114}
     E_m+CP
 = (\sum_{I\in{\cal I}^\sharp} a_I \lambda_I^{\omega} D_I \widetilde{X}_I)^{-1}.
\end{equation}
This shows  that $\det(E_m+CP)$ is a~unit of~${\cal F}$ so that $(P,C)\in\Fad$.

Next we show that four blocks of (\ref{E:H(P,C)}) are over~${\cal A}$.  The
$(2,2)$-block is the inverse of (\ref{E:08.Feb.97.224114}):
\[
     (E_m+CP)^{-1}=
    \sum_{I\in{\cal I}^\sharp} a_I \lambda_I^{\omega} D_I \widetilde{X}_I.
\]
{Similarly}, simple calculations show that other blocks are also over~${\cal A}$
as follows:%
\\
$(2,1)$-block:
\[
   C(E_n+PC)^{-1} = \sum_{I\in{\cal I}^\sharp} a_I \lambda_I^{\omega} D_I \widetilde{Y}_I,
\]
$(1,1)$-block:
\[
     (E_n+PC)^{-1}=E_n-\sum_{I\in{\cal I}^\sharp} a_I \lambda_I^{\omega} N_I\widetilde{Y}_I,
\]
$(1,2)$-block:
\[
    -P(E_m+CP)^{-1} = -\sum_{I\in{\cal I}^\sharp} a_I \lambda_I^{\omega} N_I\widetilde{X}_I.
\]

To finish the proof, we proceed to deal with the case where the matrix
$(\sum_{I\in{\cal I}^\sharp} a_I \lambda_I^{\omega} D_I \widetilde{X}_I)$ is ${\cal Z}$-singular.
To make the matrix ${\cal Z}$-nonsingular, we reconstruct the
matrices~$\widetilde{X}_I$ and~$\widetilde{Y}_I$ with an $I\in{\cal I}^\sharp$.

Since the sum of~$a_I \lambda_I^{\omega}$'s for $I\in{\cal I}^\sharp$ is equal to~$1$, by
Remark\,\ref{R:4.4}(i, iii) there exists at least one summand, say
$a_{I_0} \lambda_{I_0}^{\omega}$ with an $I_0\in{\cal I}^\sharp$, such that both $a_{I_0}$
and $\lambda_{I_0}$ belong to~${\cal A}\backslash{\cal Z}$.  Let~$R_{I_0}$ be a~parameter
matrix of ${\cal A}_{\lambda_{I_0}}^{m\times n}$.  Then the following matrix equation
holds over~${\cal A}_{\lambda_{I_0}}$:
\[%
   (\widetilde{Y}_{I_0}+R_{I_0}\widetilde{D}_{I_0})N_{I_0}
  +(\widetilde{X}_{I_0}-R_{I_0}\widetilde{N}_{I_0})D_{I_0}
  =E_m,
\]%
where $\widetilde{D}_{I_0}=\det(\lambda_{I_0}^{\omega}D_{I_0})E_n$ and $\widetilde{N}_{I_0}=\lambda_{I_0}^{\omega}N_{I_0}\adj(\lambda_{I_0}^{\omega}D_{I_0})$.
Since~$\omega$ is a~sufficiently large integer,
the following matrix equation is over~${\cal A}$:
\[%
\begin{array}{l}
   \Bigl(\lambda_{I_0}^{\omega}(\widetilde{Y}_{I_0}+R_{I_0}\widetilde{D}_{I_0})\Bigr)
   \Bigl(\lambda_{I_0}^{\omega}N_{I_0}\Bigr) \\
\hspace*{4em}
  +\Bigl(\lambda_{I_0}^{\omega}(\widetilde{X}_{I_0}-R_{I_0}\widetilde{N}_{I_0})\Bigr)
   \Bigl(\lambda_{I_0}^{\omega}D_{I_0}\Bigr)
  =\lambda_{I_0}^{2 \omega}E_m,
\end{array}
\]%
where the matrices surrounded by ``$\Bigl($'' and ``$\Bigr)$'' in the
left-hand side are over~${\cal A}$.  From the first matrix equation of
(\ref{E:09.Feb.97.151238}), $\det(D)=\det(D_I)\det(V_{\lambda_I})$ over ${\cal A}_
{\lambda_I}$ for every $I\in{\cal I}^\sharp$.  Thus by Remark\,\ref{R:4.4}(iii) the
matrix $\lambda_{I_0}^{\omega}D_{I_0}$ is ${\cal Z}$-nonsingular and so is the matrix
$\widetilde{D}_{I_0}$ ($=\det(\lambda_{I_0}^{\omega}D_{I_0})E_n$).

Consider now the following matrix equation over~${\cal A}$:
\begin{equation}\label{E:08.Feb.97.220842}
  \begin{array}{r}
    \bmatrix{
      \sum_{I\in{\cal I}^\sharp} a_I \lambda_I^{\omega} D_I\widetilde{X}_I &
      \sum_{I\in{\cal I}^\sharp} a_I \lambda_I^{\omega} D_I\widetilde{Y}_I \cr
      -a_{I_0} \lambda_{I_0}^{\omega} \det(\lambda_{I_0}^{\omega}D_{I_0})\widetilde{N}_{I_0} &
       a_{I_0} \lambda_{I_0}^{\omega} \det(\lambda_{I_0}^{\omega}D_{I_0})\widetilde{D}_{I_0}
     }
    \bmatrix{ D & O \cr
              N & E_n
     }\hspace*{4em}\mbox{}\medskip\\
    =
    \bmatrix{ D & \sum_{I\in{\cal I}^\sharp} a_I \lambda_I^{\omega} D_I\widetilde{Y}_I \cr
              O & a_{I_0} \lambda_{I_0}^{\omega} \det(\lambda_{I_0}^{\omega}D_{I_0})\widetilde{D}_{I_0} }.
  \end{array}
\end{equation}
The $(1,1)$-block of (\ref{E:08.Feb.97.220842}) can be understood in
the following way.  From the last matrix equation in
(\ref{E:09.Feb.97.151238}) we have the following matrix equation
over~${\cal A}_{\lambda_I}$:
\[
   D_I\widetilde{Y}_IN
  +D_I\widetilde{X}_ID
  =D.
\]
Considering the above equation multiplied by $a_I \lambda_I^{\omega}$ over~${\cal A}$, we
have the following equation over~${\cal A}$:
\begin{equation}\label{E:09.Feb.97.165546}
    a_I \lambda_I^{\omega} D_I\widetilde{Y}_IN
   +a_I \lambda_I^{\omega} D_I\widetilde{X}_ID
  =
    a_I \lambda_I^{\omega} D
   +a_I \lambda_I^{\omega} Z,
\end{equation}
where~$Z$ is a~matrix over~${\cal A}$ such that $\lambda_I^xZ$ is the zero matrix
for some positive integer~$x$.  Since~$\omega$ is a~large positive
integer, we can consider that the matrix $a_I \lambda_I^{\omega}Z$ in
(\ref{E:09.Feb.97.165546}) becomes the zero matrix.  Therefore, the
$(1,1)$-block of (\ref{E:08.Feb.97.220842}) holds.  Then the matrix of
the right-hand side of (\ref{E:08.Feb.97.220842}) is ${\cal Z}$-nonsingular
since both of the matrices~$D$ and 
  $a_{I_0} \lambda_{I_0}^{\omega} \det(\lambda_{I_0}^{\omega}D_{I_0})\widetilde{D}_{I_0}$ 
in the right-hand side of (\ref{E:08.Feb.97.220842}) are ${\cal Z}
$-nonsingular.  Hence the first matrix of (\ref{E:08.Feb.97.220842})
is also ${\cal Z}$-nonsingular by Remark\,\ref{R:4.4}(iii).
By Lemma\,\ref{L:4.3} and (\ref{E:08.Feb.97.220842}), there exists
a~matrix $R_{I_0}'$ of ${\cal A}^{m\times n}$ such that the following matrix is $
{\cal Z}$-nonsingular:
\begin{equation}\label{E:09.Jan.99.133634}
     \sum_{I\in{\cal I}^\sharp} a_I \lambda_I^{\omega} D_I\widetilde{X}_I
     -a_{I_0} \lambda_{I_0}^{2 \omega} D_{I_0}\adj(\lambda_{I_0}^{\omega}D_{I_0})R_{I_0}'\widetilde{N}_{I_0}.
\end{equation}
Let now 
  $R_{I_0}:=\lambda_{I_0}^{\omega}\adj(\lambda_{I_0}^{\omega}D_{I_0})R_{I_0}'$, 
  $\widetilde{X}_{I_0}:=\widetilde{X}_{I_0}-R_{I_0}\widetilde{N}_{I_0}$, and 
  $\widetilde{Y}_{I_0}:=\widetilde{Y}_{I_0}+R_{I_0}\widetilde{D}_{I_0}$.
Then the matrix $ \sum_{I\in{\cal I}^\sharp} a_I \lambda_I^{\omega}D_I\widetilde{X}_I$
becomes equal to (\ref{E:09.Jan.99.133634}) and ${\cal Z}$-nonsingular
\qquad
\end{proofof}

\begin{remark}\label{R:5.1}
From the proof above, if~(i) we can check (\ref{E:Th:3.2:1}) and
if~(ii) we can construct the right-coprime factorizations of the given
causal plant over ${\cal A}_{\lambda_I}$ for every $I\in{\cal I}^\sharp$, then we can
construct stabilizing controllers of the plant, where~$\lambda_I$ is an
element of the generalized elementary factor of the plant.  For~(i),
if we can compute, for example, the Gr\"{o}bner
basis\cite{bib:geddes92a} over~${\cal A}$ and if the generalized elementary
factors of the plant are finitely generated, (\ref{E:Th:3.2:1}) can be
checked.  For~(ii), it is already known by Lemmas\,\ref{L:4.6}
and~\ref{L:4.5} that there exist the right-coprime factorizations of
the plant over~${\cal A}_{\lambda_ I}$.

Let us give an example concerning the Gr{\"{o}}bner basis.  Consider
the generalized elementary factors of
Example\,\ref{Ex:Anantharam:3.4}.  They are expressed as
\begin{eqnarray*}
  \Lambda_{P\!I_1}&=&
       (\langle(1+2z)(1+z+z^2)(1-3z)\rangle)
      +(\langle(1+2z)(1+z+z^2)z^2\rangle)
\\&&\hspace*{5em}
      +(\langle(1+2z)(1+z+z^2)z^3\rangle),\\
  \Lambda_{P\!I_2}&=&
       (\langle(1+z)(1+2z+4z^2)(1-3z)\rangle)
      +(\langle(1+z)(1+2z+4z^2)z^2\rangle)
\\&&\hspace*{5em}
      +(\langle(1+z)(1+2z+4z^2)z^3\rangle),\\
  \Lambda_{P\!I_3}&=&
       (\langle(1+z)(1+2z)(1-3z)\rangle)
      +(\langle(1+z)(1+2z)z^2\rangle)
\\&&\hspace*{5em}
      +(\langle(1+z)(1+2z)z^3\rangle).
\end{eqnarray*}
Hence each of them has three generators and so is finitely generated.
Suppose here that we can calculate the Gr{\"{o}}bner basis over~${\cal A}$
(of Example\,\ref{Ex:Anantharam:3.4}).  Then as above the plant is
stabilizable if and only if the Gr{\"{o}}bner basis of the nine
generators contains~$1$.
\qquad\endproof
\end{remark}

In the following two examples we follow the proof of
Theorem\,\ref{Th:3.3}.  In the first one, we construct a~stabilizing
controller with part~(c).  In the other example, we follow part~(b).  On the other hand we do not follow part~(a) since it
can be followed easily with part~(a) of~(\romref{P:4.1:i}) in the
proof of Proposition\,\ref{P:4.1}.

\begin{example}\label{Ex:Anantharam:5.3}(Continued)  We continue 
Example\,\ref{Ex:Anantharam:3.4} (and \ref{Ex:12.Nov.99.111204}) and
construct a~stabilizing controller of the plant as in
``(\romref{Th:3.3:iii})$\rightarrow$(\romref{Th:3.3:i})'' of the proof above.
Let the notation be as in Examples\,\ref{Ex:Anantharam:3.4}
and~\ref{Ex:12.Nov.99.111204}.

Since, in this example, $\Lambda_{P\!I_1}+\Lambda_{P\!I_2}={\cal A}$ holds, ${\cal I}^\sharp=
\{ I_1,I_2\}$.  For $I_1\in{\cal I}^\sharp$, the matrices~$N_{I_1}$, $D_{I_1}$,
$\widetilde{X}_{I_1}$ and $\widetilde{Y}_{I_1}$ of
(\ref{E:09.Feb.97.151238}) over~${\cal A}_{\lambda_{I_1}}$ have been calculated
as (\ref{E:02.Nov.99.150110:1}) and (\ref{E:02.Nov.99.150110:2}).  For
$I_2\in{\cal I}^\sharp$, the matrices~$N_{I_2}$, $D_{I_2}$, $\widetilde{X}_{I_2}$
and $\widetilde{Y}_{I_2}$ of (\ref{E:09.Feb.97.151238}) over~${\cal A}_{\lambda_
{I_2}}$ can be calculated analogously as follows:
\begin{eqnarray*}
  N_{I_2}&=&\bmatrix{\lambda_{I_2}^{-1}\alpha_{I_2}\langle(1+2z)(1-3z+z^2)(1+z+z^2)\rangle \cr
                     1},\\
  D_{I_2}&=&\bmatrix{\lambda_{I_2}^{-1}\alpha_{I_2}\langle(1+z)(1+2z)(1-3z+z^2)\rangle},\\
  \widetilde{Y}_{I_2}&=& \bmatrix{ 0 & 1 }, 
  \widetilde{X}_{I_2} =  \bmatrix{ 0 }.
\end{eqnarray*}

Then the following matrices are over~${\cal A}$:
\begin{eqnarray*}
 \lambda_{I_1}D_{I_1}\widetilde{X}_{I_1}&=&
        \bmatrix{ 0 },{\ \ }
 \lambda_{I_1}D_{I_1}\widetilde{Y}_{I_1} =
        \bmatrix{ \alpha_{I_1}\langle(1+z)(1+2z)(1-3z)\rangle & 0 },
\\
 \lambda_{I_2}D_{I_2}\widetilde{X}_{I_2}&=&
        \bmatrix{ 0 },{\ \ }
 \lambda_{I_2}D_{I_2}\widetilde{Y}_{I_2}=
        \bmatrix{ 0 & \alpha_{I_2}\langle(1+z)(1+2z)(1-3z+z^2)\rangle }.
\end{eqnarray*}
Hence in this example, we can let $\omega=1$ as a~sufficiently large
integer and $a_I=1$ for all $I\in{\cal I}^\sharp$ (since $ \sum_{I\in{\cal I}^\sharp}\lambda_I^{\omega}%
=1$).

Note here that the matrix $\lambda_{I_1}D_{I_1}\widetilde{X}_{I_1} +\lambda_
{I_2}D_{I_2}\widetilde{X}_{I_2}$ is ${\cal Z}$-singular.  Hence we should
reconstruct the matrices $\widetilde{Y}_{I_i}$ and
$\widetilde{X}_{I_i}$ with~$i$ being either~$1$ or~$2$ as in the proof
of Theorem\,\ref{Th:3.3}.  Since, in this example, both $\lambda_{I_1}$ and
$\lambda_{I_2}$ are nonzerodivisors, we can choose each of~$1$ and~$2$.
This example proceeds by reconstructing the
matrices~$\widetilde{Y}_{I_1}$ and~$\widetilde{X}_{I_1}$, which means
that~$I_1$ is used as~$I_0$ in the proof of Theorem\,\ref{Th:3.3}.
The actual reconstruction is done by following the proof of
Lemma\,\ref{L:4.3}.

Consider the first matrix of (\ref{E:08.Feb.97.220842}).  Recall that
  $\widetilde{N}_{I_0}=\lambda_{I_0}^{\omega}N_{I_0}\adj(\lambda_{I_0}^{\omega}D_{I_0})$ and
  $\widetilde{D}_{I_0}=\det(\lambda_{I_0}^{\omega}D_{I_0})E_n$.
In this example, they are given as
\begin{eqnarray*}
  \widetilde{N}_{I_1}&=&(\widetilde{N}_{I_0}=)
                        \bmatrix{ \lambda_{I_1}  \cr
                                  \alpha_{I_1}\langle(1+z)(1-3z)(1+2z+4z^2)\rangle},\\
  \widetilde{D}_{I_1}&=&(\widetilde{D}_{I_0}=)
                        \alpha_{I_1}\langle(1+z)(1+2z)(1-3z)\rangle E_2.
\end{eqnarray*}
One can check that the first matrix of (\ref{E:08.Feb.97.220842})
is ${\cal Z}$-nonsingular.  Then we construct a~matrix~$R_{I_0}'$ of ${\cal A}^{1
\times 2}$ such that (\ref{E:09.Jan.99.133634}) is ${\cal Z}$-nonsingular.  To do
so, we follow temporarily the proof of the Lemma\,\ref{L:4.3}.

Consider the first matrix of (\ref{E:08.Feb.97.220842}) as the matrix
of (\ref{E:L:4.3:1}), that is,
\begin{eqnarray*}
 A &=& \sum_{I\in{\cal I}^\sharp} a_I \lambda_I^{\omega} D_I\widetilde{X}_I=\bmatrix{0},\\
 B &=& -a_{I_0} \lambda_{I_0}^{\omega} \det(\lambda_{I_0}^{\omega}D_{I_0})\widetilde{N}_{I_0}\\
   &=& -\alpha_{I_1}\lambda_{I_1}\langle(1+z)(1+2z)(1-3z)\rangle
        \bmatrix{ \lambda_{I_1} \cr
                  \alpha_{I_1}\langle(1+z)(1-3z)(1+2z+4z^2)\rangle}.
\end{eqnarray*}
Then we choose a~full-size~$a$ minor of $\bmatrix{A^t & B^t}^t$ having
as few rows from~$B$ as possible.  In this example, we can choose both
entries in~$B$.  Here we choose the $(1,1)$-entry of~$B$, so that
\begin{equation}\label{E:01.Nov.99.182529}
   a=-\alpha_{I_1}\lambda_{I_1}^2\langle(1+z)(1+2z)(1-3z)\rangle.
\end{equation}
Thus we have $k=1$, $i_1=1$ and $j_1=1$, where the notations~$k$,
$i_1,\ldots, i_k$, and $j_1,\ldots, j_k$ are as in the proof of
Lemma\,\ref{L:4.3}.  Hence~$R$ in the proof is given as $R=\bmatrix{1
& 0}$.  We can confirm that $A+RB=\bmatrix{a}$ which is ${\cal Z}
$-nonsingular by observing that every factor of the right-hand side of
(\ref{E:01.Nov.99.182529}) has a~nonzero constant term.

From here on we proceed with following again the proof of
Theorem\,\ref{Th:3.3}.  The notation~$R$ used above corresponds to the
notation $R_{I_0}'$ in the proof of Theorem\,\ref{Th:3.3} (that is,
$R_{I_0}'=\bmatrix{1 & 0}$).  The matrix $R_{I_1}$ is given as
follows:
\[
  R_{I_1} =
 (R_{I_0} =)\lambda_{I_1}^{\omega}\adj(\lambda_{I_1}^{\omega}D_{I_1})R_{I_1}'
          = \lambda_{I_1}\bmatrix{1 & 0}.
\]
Then new $\widetilde{X}_{I_1}$ and $\widetilde{Y}_{I_1}$ are given as
follows:
\begin{eqnarray*}
  \widetilde{X}_{I_1}&:=&\widetilde{X}_{I_1}-R_{I_1}\widetilde{N}_{I_1}
                       = \bmatrix{ -\lambda_{I_1}^2 },\\
  \widetilde{Y}_{I_1}&:=&\widetilde{Y}_{I_1}+R_{I_1}\widetilde{D}_{I_1}
                       = \bmatrix{1+\alpha_{I_1}\lambda_{I_1}\langle(1+z)(1+2z)(1-3z)\rangle&0}.
\end{eqnarray*}
Therefore a~stabilizing controller~$C$ of the form
(\ref{E:08.Feb.97.215931}) is obtained as
\begin{eqnarray*}
 C&=& (\lambda_{I_1}D_{I_1}\widetilde{X}_{I_1}
     + \lambda_{I_2}D_{I_2}\widetilde{X}_{I_2})^{-1}
      (\lambda_{I_1}D_{I_1}\widetilde{Y}_{I_1}
     + \lambda_{I_2}D_{I_2}\widetilde{Y}_{I_2})
\\
  &=&
     \frac{-1}{\alpha_{I_1}\lambda_{I_1}^2\langle(1+z)(1+2z)(1-3z)\rangle}
\\
  & &
~~~
     \bmatrix{\alpha_{I_1}\langle(1+z)(1+2z)(1-3z)\rangle(1+\alpha_{I_1}\lambda_{I_1}\langle(1+z)(1+2z)(1-3z)\rangle) \cr
              \alpha_{I_2}\langle(1+z)(1+2z)(1-3z+z^2)\rangle}^t.
\end{eqnarray*}

The matrix $H(P,C)$ with the stabilizing controller~$C$ above over~${\cal A}
$ is expressed as follows:
\[
   H(P,C)=\bmatrix{ h_{11} & h_{12} & h_{13} \cr
                    h_{21} & h_{22} & h_{23} \cr
                    h_{31} & h_{32} & h_{33}},
\]
where 
\begin{eqnarray*}
  h_{11}&=& -\alpha_{I_1}\lambda_{I_1}^2\langle(1+z)(1+2z)(1-3z)\rangle\\
        & & \hspace*{0.15\textwidth} +\alpha_{I_2}\langle(1+z)(1-3z+z^2)(1+2z+4z^2)\rangle,\\
  h_{12}&=&-\alpha_{I_2}\langle(1+2z)(1+z+z^2)(1-3z+z^2)\rangle,\\
  h_{13}&=& \lambda_{I_1}^3,\\
  h_{21}&=&-\alpha_{I_1}\langle(1+z)(1-3z)(1+2z+4z^2)\rangle\\
        & & \hspace*{0.15\textwidth}  (1+\lambda_{I_1}\alpha_{I_1}\langle(1+z)(1+2z)(1-3z)\rangle),\\
  h_{22}&=& \alpha_{I_1}(\langle(1+2z)(1-3z)(1+z+z^2)\rangle
                         (1+\alpha_{I_1}\lambda_{I_1}\langle(1+z)(1+2z)(1-3z)\rangle)\\
        & & \hspace*{0.15\textwidth} -\lambda_{I_1}^2\langle(1+z)(1+2z)(1-3z)\rangle),\\
  h_{23}&=& \alpha_{I_1}\lambda_{I_1}^2\langle(1+z)(1-3z)(1+2z+4z^2)\rangle,\\
  h_{31}&=& \alpha_{I_1}\langle(1+z)(1+2z)(1-3z)\rangle(1+\alpha_{I_1}\lambda_{I_1}\langle(1+z)(1+2z)(1-3z)\rangle),\\
  h_{32}&=& \alpha_{I_2}\langle(1+z)(1+2z)(1-3z+z^2)\rangle,\\
  h_{33}&=&-\alpha_{I_1}\lambda_{I_1}^2\langle(1+z)(1+2z)(1-3z)\rangle.
\end{eqnarray*}

Before finishing this example, let us show that the ${\cal A}$-module~${\cal T}_P$
is \emph{not} free.  We show it by contradiction.  Suppose that~${\cal T}_P$
is free.  Then the ${\cal A}$-module $M_r(T)$ is also free.  Since the
matrix~$D$, a~part of~$T$, is nonsingular, the rank of $M_r(T)$
is~$m$.  Let~$V$ be a~matrix in $({\cal A})_m$ whose rows are~$m$ distinct
generators of $M_r(T)$.  As in (\ref{E:09.Feb.97.151238}), we have
matrices $\widetilde{Y}$, $\widetilde{X}$, $N'$, $D'$ over~${\cal A}$ such
that
\[%
  \bmatrix{N^t & D^t}^t=\bmatrix{N'^t & D'^t}^t V,{\ \ }
  \bmatrix{\widetilde{Y} & \widetilde{X}}\bmatrix{N^t & D^t}^t=V,{\ \ }
   \widetilde{Y}'N'+\widetilde{X}'D'=E_1.
\]%
However the last matrix equation is inconsistent with the fact that
the plant~$P$ does not have coprime factorization.  Therefore~${\cal T}_P$
is not free.  Nevertheless we note that~${\cal T}_P$ is projective by
Theorem\,\ref{Th:3.3}.
\qquad\endproof
\end{example}

\begin{example}(Continued)
Let us follow part~(b) in the proof of Theorem\,\ref{Th:3.3}.
Suppose that (\romref{Th:3.3:i}) of Theorem\,\ref{Th:3.3} holds, that
is, the modules~${\cal T}_P$ and ${\cal W}_ P$ are projective.

Consider again the plant~$P$ of (\ref{E:3.2}).  Let $F_1=\{\lambda_{I_1},\lambda_
{I_2}\}$, where $\lambda_{I_1}$ and $\lambda_{I_2}$ are given as in
(\ref{E:12.Nov.99.131258}).  Then we have known that $\Sigma_{f\in F_1}f=1$
and that there exists a~right-coprime factorization of the plant
over~${\cal A}_f$ for every $f\in F_1$.  By Lemma\,\ref{L:4.2}, the transposed
plant~$P^t$ is stabilizable.  We can construct its stabilizing
controller by analogy to Example~\ref{Ex:Anantharam:5.3}.  Further we
see that for both $\lambda_{I_1}$ and $\lambda_{I_2}$, the transposed
plant~$P^t$ has right-coprime factorizations over ${\cal A}_{\lambda_{I_1}}$ and
${\cal A}_{\lambda_{I_2}}$; that is, $P$ has left-coprime factorizations over ${\cal A}_
{\lambda_{I_1}}$ and ${\cal A}_{\lambda_{I_2}}$.  Thus let $F_2=\{\lambda_{I_1},\lambda_{I_2}\}
$.  For $\lambda_{I_1}\in F_2$, we have the matrices $\widetilde{N}_{I_1}$
$\widetilde{D}_{I_1}$, $Y_{I_1}$, $X_{I_1}$ over ${\cal A}_{\lambda_{I_1}}$ such
that $\widetilde{N}_{I_1}Y_{I_1}+\widetilde{N}_{I_1}X_{I_1}=E_2$ and
\begin{eqnarray*}
  \widetilde{N}_{I_1}&=& \textstyle [{ 1 \atop 0 }],{\ \ }
  Y_{I_1} = \bmatrix{ 1 & 0 },{\ \ }
  X_{I_1} = \textstyle [{ 0~0 \atop 0~1}],\\
  \widetilde{D}_{I_1}&=&
      \bmatrix{ \lambda_{I_1}^{-1}\alpha_{I_1}\langle(1+z)(1+2z)(1-3z)\rangle & 0 \cr
                \lambda_{I_1}^{-1}\alpha_{I_1}\langle(1+z)(1-3z)(1+2z+4z^2)\rangle & 1}.
\end{eqnarray*}
On the other hand, for $\lambda_{I_2}\in F_2$, we have the matrices
$\widetilde{N}_{I_2}$ $\widetilde{D}_{I_2}$, $Y_{I_2}$, $X_{I_2}$ over
${\cal A}_{\lambda_{I_2}}$ such that
$\widetilde{N}_{I_2}Y_{I_2}+\widetilde{N}_{I_2}X_{I_2}=E_2$ and
\begin{eqnarray*}
  \widetilde{N}_{I_2}&=& \textstyle [{ 1 \atop 0 }],{\ \ }
  Y_{I_2} = \bmatrix{ 1 & 0 },{\ \ }
  X_{I_2} = \textstyle [{ 0~1 \atop 0~0}],\\
  \widetilde{D}_{I_2}&=&
      \bmatrix{ 0 &  \lambda_{I_2}^{-1}\alpha_{I_2}\langle(1+z)(1+2z)(1-3z+z^2)\rangle \cr
                1 & -\lambda_{I_2}^{-1}\alpha_{I_2}\langle(1+2z)(1+z+z^2)(1-3z+z^2)\rangle}.
\end{eqnarray*}
Now we let $F=\{\lambda_{I_1}^2,\lambda_{I_1}\lambda_{I_2},\lambda_{I_2}^2\}$ ($=\{ f_1f_2\,|\,
f_1\in F_1,f_2\in F_2\}$).  Then~$F$ still generates~${\cal A}$ since $\lambda_
{I_1}^2+2 \lambda_{I_1}\lambda_{I_2}+\lambda_{I_2}^2=1$.

In the following we consider the case $f=\lambda_{I_1}^2$.  Then using the
matrix~$K$ of (\ref{E:11.Nov.99.213238}), we have
(\ref{E:11.Nov.99.213521}) with $\nu=1$, $\Kf=K$ and $V_f=\lambda_{I_1}\Delta_
{I_1}T$.

Then the ideal $I_{m{\cal A}}(f^{\nu}\Kf)$ is generated by
\[
  \lambda_{I_1}^3,{\ \ }
  \alpha_{I_1}\lambda_{I_1}^2\langle(1+z)(1-3z)(1+2z+4z^2)\rangle,{\ \ }
  \alpha_{I_1}\lambda_{I_1}^2\langle(1+z)(1+2z)(1-3z)\rangle.
\]
Thus since each of them is in $\Lambda_{P\!I_1}$,
(\ref{E:14.Feb.97.113224}) holds.  Further we can observe that for any
integer~$\xi$ greater than~$1$, (\ref{E:08.Feb.97.134530}) holds since
$\lambda_{I_1}^3\in I_{m{\cal A}}(f^{\nu}\Kf)$.

For the other cases $f=\lambda_{I_1}\lambda_{I_2}$ and $f=\lambda_{I_2}^2$, we can
follow the relations of (\ref{E:14.Feb.97.113224}) and
(\ref{E:08.Feb.97.134530}) analogously.  Details are left to
interested readers.  \qquad\endproof
\end{example}

\begin{remark}
Since Anantharam's example in~\cite{bib:anantharam85a} is artificial,
we do not present here the construction of a~stabilizing controller.
However we can construct it as part~(c) in the proof of
Theorem\,\ref{Th:3.3} (Since Anantharam in~\cite{bib:anantharam85a}
did not consider the causality, we let ${\cal Z}=\{ 0\}$ so that ${\cal P}={\cal F}$).
\end{remark}

\section{Causality of Stabilizing Controllers}
\label{S:CSC}
In this section, we present two facts: (i) for a~stabilizable causal
plant, there exists at least one stabilizing causal controller and~(ii) the stabilizing controller of the strictly causal plant is
causal, which inherits~Theorem\,4.1 in~\S\,III
of~\cite[p.888]{bib:vidyasagar82a} and Proposition\,1
of~\cite{bib:sule94a}.

\begin{proposition}\label{P:6.1}
For every stabilizable causal plant, there exists at least one
stabilizing causal controller of the plant.
\end{proposition}
\begin{proof}
In the construction of the stabilizing controller in part~(c) of
the proof of Theorem\,\ref{Th:3.3}, the denominator matrix of
(\ref{E:08.Feb.97.215931}) is ${\cal Z}$-nonsingular.  Suppose that the
obtained stabilizing controller is expressed as
$\widetilde{B}^{-1}\widetilde{A}$ with the matrices~$\widetilde{A}$
and~$\widetilde{B}$ over~${\cal A}$ such that~$\widetilde{B}$ is ${\cal Z}
$-nonsingular.  Then since the relation
  $\widetilde{B}^{-1}\widetilde{A}
        =(\det(\widetilde{B})E_m)^{-1}(\adj(\widetilde{B})\widetilde{A})$ 
holds, every entry of $\widetilde{B}^{-1}\widetilde{A}$ is causal.
\qquad
\end{proof}

\begin{proposition}\label{P:6.2}
For every stabilizable strictly causal plant, all stabilizing
controllers of the plant must be causal.
\end{proposition}
\begin{proof}
Suppose that the plant~$P$ is stabilizable and strictly causal.
Suppose further that~$C$ is a~stabilizing controller of~$P$.  We
employ the notation from part~(c) of the proof of
Theorem\,\ref{Th:3.3}.  Thus, $a_{I_0},\lambda_{I_0}\in{\cal A}\backslash{\cal Z}$ and
 $ \widetilde{Y}_{I_0}N_{I_0}
  +\widetilde{X}_{I_0}D_{I_0}
  =E_m$
with $P=N_{I_0}D_{I_0}^{-1}\in{\cal F}({\cal A}_{\lambda_{I_0}})$ from
(\ref{E:09.Feb.97.151238}).  Let
  ${\cal Z}_{\lambda_{I_0}}=\{ z/1\cdot u\in{\cal A}_{\lambda_{I_0}}\,|\, z\in{\cal Z}, 
                          \mbox{$u$ is a~unit of ${\cal A}_{\lambda_{I_0}}$}\}$.
Then this ${\cal Z}_{\lambda_{I_0}}$ is again a~principal ideal of ${\cal A}_{\lambda_
{I_0}}$.

Observe here that Lemma\,8.3.2 of~\cite{bib:vidyasagar85a} and its
proof hold even over a~general commutative ring.  According to its
proof, there exist matrices $\widetilde{A}$ and $\widetilde{B}$ over $
{\cal A}_{\lambda_{I_0}}$ such that $C=\widetilde{B}^{-1}\widetilde{A}$ and
 $ \widetilde{A}N_{I_0}
  +\widetilde{B}D_{I_0}
  =E_m$
($\widetilde{A}$ and $\widetilde{B}$ correspond to~$T$ and~$S$,
respectively, in the proof of Lemma\,8.3.2
of~\cite{bib:vidyasagar85a}).  Observe also that every entry of
$N_{I_0}$ is in ${\cal Z}_{\lambda_{I_0}}$.  Thus reviewing the proof of
Lemma\,3.5 of~\cite{bib:vidyasagar82a}, in which the calligraphic~$H$
and~$K$ in~\cite{bib:vidyasagar82a} correspond to~${\cal A}_{\lambda_{I_0}}$
and~${\cal Z}_{\lambda_{I_0}}$, respectively, we have $\det\widetilde{B}\in{\cal A}_{\lambda_
{I_0}}\backslash{\cal Z}_{\lambda_{I_0}}$.  This implies that
$\widetilde{B}^{-1}\widetilde{A}\in{\cal P}^{m\times n}$ by noting that $\lambda_{I_0}
\in {\cal A}\backslash{\cal Z}$.  Thus~$C$ is causal.  \qquad
\end{proof}

\section{Further Work}
In this paper we have presented criteria for feedback stabilizability.
We have also presented a~construction of a~stabilizing controller to
which Sule's method cannot be applied.  Recently the first
author\cite{bib:mori99ds} has developed a~parameterization of
stabilizing controllers, which is based on the results of this paper
and which does not require coprime factorizability.  This can be
applied to models to which Youla-Ku\v{c}era parameterization
cannot be applied.

\end{document}